\documentclass{amsart}
\usepackage[english]{babel}

\usepackage{amsmath, amssymb, amsfonts, amsthm, mathtools, enumerate}

\usepackage[
colorlinks=true,
linkcolor=blue!70!black,
citecolor=blue!70!black,
urlcolor=blue!70!black
]{hyperref}

\usepackage{pgfplots}
\pgfplotsset{compat=1.18}
\usepackage{tikz}
\usepackage{graphics}

\newtheorem*{theorem*}{Main Theorem}
\newtheorem{theorem}{Theorem}[section]
\newtheorem{definition}[theorem]{Definition}
\newtheorem{lemma}[theorem]{Lemma}
\newtheorem{prop}[theorem]{Proposition}

\newtheorem{remark}[theorem]{Remark}
\newtheorem{example}[theorem]{Example}

\newtheorem{innercustomthm}{}
\newenvironment{customthm}[1]
  {\renewcommand\theinnercustomthm{#1}\innercustomthm}
  {\endinnercustomthm}

\newcommand{\BR}{\mathbb{R}}

\newcommand{\BN}{\mathbb{N}}

\newcommand{\R}{\BR}
\newcommand{\e}{\varepsilon}
\newcommand{\id}{\mathrm{id}}

\newcommand{\mcP}{\mathcal P}

\newcommand{\dist}{\mathop{\mathrm{dist}}}

\newcommand{\abs}[1]{\left|#1\right|}
\newcommand{\Abs}[1]{\Big|#1\Big|}

\newcommand{\der}[2][\tau]{\frac{\partial #2}{\partial #1}}
\newcommand{\derder}[2][\tau]{\frac{\partial^2 #2}{\partial #1^2}}
\newcommand{\osc}{\mathop{\mathbf{osc}}}
\newcommand{\const}{\mathcal{O}(1)}
\newcommand{\ind}[1]{\mathbf{1}_{#1}}

\newcommand{\restr}[1]{\big|_{#1}}

\newcommand{\TV}{\mathop{\mathrm{TV}}\nolimits}

\newcommand{\Loneloc}{L^1_{\mathrm{loc}}}
\newcommand{\norm}[1]{\left\|#1\right\|}
\newcommand{\BV}{\mathop{\mathrm{BV}}}

\newcommand{\RcOp}{ \mathbf{Cut}}
\newcommand{\Rc}[2]{\RcOp_{#1, #2}}
\newcommand{\Cr}[1]{\mathbf{Crop}_{#1}}

\newcommand{\conv}{\mathop{\mathbf{conv}}}

\newcommand{\Height}{\mathbf{H}}
\newcommand{\cost}{\mathop{\mathbf{cost}}\nolimits}
\newcommand{\bleft}{ b^{\rm left}}
\newcommand{\bright}{b^{\rm right}}
\newcommand{\stress}{\mathop{\mathbf{stress}}}
\newcommand{\bridges}{\mathop{\mathbf{Br}}\nolimits}
\newcommand{\PL}{\mathrm{PL}}

\newcommand{\psishape}{\mathop{%
				\tikz[baseline =-0.1ex, x=2ex, y=2ex, yscale=1.7]{
					\draw[line width = 0.1ex, line cap=round, rounded corners = 0.025ex]
					plot[domain=0.01:0.5,samples=10]
					(\x, {\x - \x*ln(2*\x)}) -- plot[domain=0.5:0.99,samples=10]
					(\x, {1-\x - (1-\x)*ln(2-2*\x)})
					-- (1,0) -- (0,0) -- cycle;
			}\mkern 1mu}\nolimits}

\newcommand{\discrpar}{\square}

\newcommand{\mcF}{\mathcal F}
\newcommand{\mcD}{\mathcal D}

\newcommand{\Qcost}{\mathcal{Q}^{\mathrm{cost}}}

\author{Egor Baitenov}
\author{Stefano Modena}
\address{GSSI, I-67100 L'Aquila, Italy}
\email{egor.baitenov@gssi.it, stefano.modena@gssi.it}
\date{\today}

\begin{document}
\title[Local-in-time quadratic functional]{A local-in-time quadratic functional  \\ for scalar conservation laws}

\begin{abstract}
   We construct a Glimm-type functional for scalar conservation laws, which is quadratic (w.r.t. the total variation of the initial datum), it controls the change of wavefront speed across interactions,  and, differently than functionals with similar properties already present in the literature, it is purely local in time, meaning that it depends only on the present profile of the (approximate) solution (and not on its past or future evolution). 
\end{abstract}

\maketitle

\section{Introduction}

Consider the  $n\times n$ system of conservation laws
\begin{equation}
    \label{eq:syst-cons-law}
    u_t + f(u)_x  = 0, 
\end{equation}
for the unknown  $u\colon[0, + \infty) \times \R \to \R^n$
with initial datum
\begin{equation}
    \label{eq:initial-datum}
    u|_{t=0} = \bar u.
\end{equation}
Here the flux function
\begin{equation}
    f:  \overline{\Omega} \subseteq \R^n \to \R^n
\end{equation}
is smooth, with uniformly bounded second derivative, $\Omega \subseteq \R^n$ is open, and $f$ is only assumed to be strictly hyperbolic.\footnote{\emph{Strict hyperbolicity} means that the matrix $\nabla f(u)$ has $n$ distinct real eigenvalues $\lambda_1(u)<\ldots<\lambda_n(u)$ for each $u \in \Omega$.} No convexity-like assumptions\footnote{E.g.  genuine nonlinearity (GNL) or linear degeneracy (LD).} on $f$ will be made throughout. By the theory of Bianchini--Bressan~\cite{BianchiniBressan2005}, there is a subset $\mcD \subseteq \Loneloc(\R; \R^n)$, containing all $\BV$ maps with sufficiently small total variation and left state $u(-\infty)$ belonging to a prescribed compact $K \subset \Omega$, and a Lipschitz semigroup
\begin{equation*}
    S\colon[0,  +\infty) \times \mcD \to \mcD, \quad (t,\bar u) \mapsto S_t \bar u,
\end{equation*}
such that for every $\bar u \in \mcD$,  $u(t,x) := (S_t \bar u)(x)$ is the vanishing viscosity solution to~\eqref{eq:syst-cons-law} with initial datum $\eqref{eq:initial-datum}$. 

The main goal of this paper (for the scalar conservation law, $n=1$) and the companion paper \cite{BaitMod2} (for systems of conservation laws, $n \geqslant 1$) is to prove the following Theorem, stated for the moment in a quite informal  way (the rigorous statements corresponding to points~\eqref{thm-main:p1}--\eqref{thm-main:p4} below will be presented later in Theorems~\ref{thm:quadr-bd}--\ref{thm:const-if-no-interaction} and~\ref{thm:further-properties}, respectively). 
\begin{customthm}{Main Theorem}
\label{thm-main}
    There is a \emph{local-in-time} functional
    \begin{equation*}
        \Upsilon \colon \mathcal{D} \to \R_{\geqslant 0},
    \end{equation*}
    such that
    \begin{enumerate}[(i)]
        \item\label{thm-main:p1} {\underline{Quadratic bound}:}
        $\Upsilon(u) = \const \TV(u)^2$ 
        (i.e. $\Upsilon$ is quadratic w.r.t. the  total variation of its argument),\footnote{The (standard) notations $A = \const B$, $A \leqslant \const B$ means that $A \leqslant CB$ for some constant $C>0$ depending only the flux function $f$ in~\eqref{eq:syst-cons-law}.}
        \item\label{thm-main:p2} {\underline{Interaction estimate}:}
        in each \emph{interaction between wavefronts} $\Upsilon$ decreases, and its decrease bounds  the \emph{change of speed} of the wavefronts involved in the interaction,
        \item\label{thm-main:p3}  {\underline{Persistency}:}  $\Upsilon(S_t u) = \Upsilon(u)$, if $u$ is piecewise constant (with uniform spatial separation of jumps) and $t$ is small enough.
    \end{enumerate}
        Moreover,
    \begin{enumerate}[(i)]
    \setcounter{enumi}{3}
        \item\label{thm-main:p4}
        $\Upsilon$ satisfies additional properties: approximation properties, monotonicity a\-cross Riemann cuts, Lyapunov property, lower semicontinuity,  stability w.r.t. perturbation of state curve,   invariance w.r.t. reparametrization of state curve.
    \end{enumerate}
\end{customthm}

The notions of \emph{interaction between wavefronts} and \emph{change of speed} will be informally introduced in Section~\ref{ss:basic-pbl} and rigorously in Section~\ref{ss:non-neg-quadr-bound}. 

Quadratic functionals for systems of conservation laws, satisfying some version of Points~\eqref{thm-main:p1},~\eqref{thm-main:p2} and~\eqref{thm-main:p3} above, were already present in the literature (\cite{AM2011, BM-scalar, BM-CaseStudy, BM-CMP, M-Forward}). The main novelty of our functional consists in the fact that it is \emph{local-in-time}. We will comment below in Section~\ref{ss:relevance} on the meaning of the expression ``\emph{local-in-time}'' and on its relevance. 

The additional properties listed in Point~\eqref{thm-main:p4} are as well consequence (though nontrivial ones) of the locality in time of~$\Upsilon$.

We develop in this paper the theory for the scalar case (i.e. $n=1$),  defining the \emph{local-in-time} functional $\Upsilon$ and  showing that it satisfies Points~\eqref{thm-main:p1}--\eqref{thm-main:p3}
above. In the companion paper~\cite{BaitMod2}, we develop the general theory for arbitrary $n$, extending the definition of $\Upsilon$ to the case of systems and proving all Points~\eqref{thm-main:p1}--\eqref{thm-main:p3} as well as~\eqref{thm-main:p4} above.\footnote{Notice that the properties listed in~\eqref{thm-main:p4} hold as well in the scalar case $n=1$ (see Section \ref{ss:additional-properties}), but we  prove them directly in the general system case $n \geqslant 1$ in \cite{BaitMod2}.}

We chose to isolate the scalar analysis from the system one because the main ideas of our construction are purely \emph{scalar}.  Extension to systems is essentially a combination of such scalar ideas with well known interaction estimates. We thus hope that the present, scalar paper conveys, in the clearest possible way, the core of the construction, which otherwise (in the case of systems) might remain hidden behind technicalities. 

In the rest of this introduction, we discuss the relevance of interaction estimates and interaction functionals in the theory of hyperbolic systems of conservation laws, with a particular emphasis on the importance of \emph{local-in-time functionals}. Even though the rest of the paper will deal only with the scalar case, we keep the discussion in this introduction at the general system level.  

\subsection{The basic problem: estimating interactions in approximate solutions}
\label{ss:basic-pbl}

In the analysis of solutions to hyperbolic conservation laws, one of the central problems is to understand how the \emph{strength} or the \emph{speed} of waves are affected by the interactions between different wavefronts, due to the nonlinearity of the flux~$f$. 

More precisely, let us assume that we are given an approximate solution $u^\e$,  constructed e.g. by means of the wavefront tracking or the Glimm scheme, with~$\e$ denoting the approximation parameter. Such approximate solution consists of a (locally) finite number of wavefronts which interact on a discrete set of points $\{P_j\}_{j \in \BN}$, $P_j = (t_j, x_j)$, in the $(t,x)$-plane. We will use the term \emph{wavefront} in a~not rigorous way, meaning a solution (or an approximate solution) to a Riemann problem. In each interaction point $P_j$, the wavefronts involved in the interaction undergo a change of their strength and of their speed, which we denote here by $\Delta s(P_j)$ (for the change of strength) and by $\Delta \sigma(P_j)$ (for the change of speed). A~precise definition (at least in the scalar case) of these quantities will follow shortly. 
We assume also, for the sake of simplicity, that $t_j \neq t_{j'}$ if $j \neq j'$, i.e. at most one interaction takes place at each time. 

The problem of ``estimating  interactions'' in approximate solutions consists essentially in proving  bounds of the form
\begin{subequations}
\label{eq:intro-interaction-estimate}
\begin{equation}
\label{eq:intro-strength}
    \sum_j \Delta s(P_j) \leqslant C(f, \TV(\bar u)), 
\end{equation}
\begin{equation}
\label{eq:intro-speed}
    \sum_j \Delta \sigma(P_j) \leqslant C(f, \TV(\bar u)),
\end{equation}
\end{subequations}
where $C(f, \TV(\bar u))$ denotes a constant depending on the flux $f$, on the total variation of the initial datum $\bar u$, but not on the approximation parameter $\e$.

The first estimate (the one on the change of strength) is of fundamental importance for proving e.g. a (uniform in $\e$) bound on the total variation of approximate solutions, and thus to get $L^1$ compactness. The second estimate (the one on the change of speed) is relevant, on the one side, in order to get consistency of the deterministic version of the Glimm scheme~\cite{Liu1977, BM-Glimm}; on the other side, it is also useful in order to get additional compactness properties, and we will comment more extensively in Section~\ref{ss:relevance}  about this.

\begin{figure}[h]
		\centering
		\begin{minipage}{0.45\textwidth}
			\centering
			\includegraphics{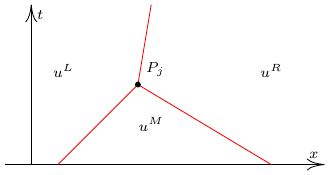}
            \\
			(a) 
            \\
            \includegraphics{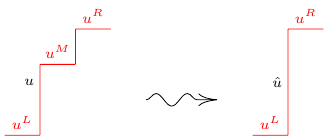}
            \\
			(c) 
		\end{minipage}
        \hfill
        \begin{minipage}{0.45\textwidth}
			\centering
			\includegraphics{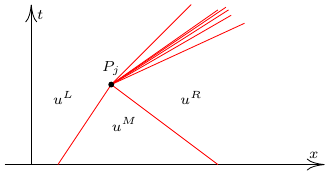}
            \\
			(b) 
            \\ \includegraphics{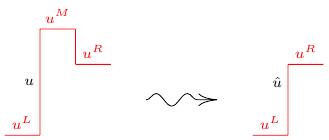}
            \\
			(d) 
		\end{minipage}
		\caption{}
		\label{fig:interaction-t-x}
	\end{figure}

We give now the  definitions of \emph{change of strength} $\Delta s(P_j)$ and \emph{change of speed} $\Delta \sigma(P_j)$ across an interaction, in the scalar case $n=1$ (in the case of systems the definition is more technical). Fix a point $P_j$ on the $(t,x)$-plane where an interaction between different wavefronts take place, see Fig.~\ref{fig:interaction-t-x}(a)--\ref{fig:interaction-t-x}(b). W.l.o.g. we can assume that before the interaction there are two wavefronts connecting respectively $u^L$ to $u^M$ and $u^M$ to $u^R$, with $u^L, u^M, u^R \in \R$,  whereas after the interaction there is a single wavefront connecting $u^L$ to $u^R$. 

There are two basic situations: the two interacting wavefronts have the same sign (w.l.o.g. $u^L \leqslant u^M \leqslant u^R$) or they have opposite sign (w.l.o.g.  $u^L \leqslant u^R \leqslant u^M$), see Fig.~\ref{fig:interaction-t-x}(c)--\ref{fig:interaction-t-x}(d). 

The \emph{change of strength} across the interaction is the quantity
\begin{equation*}
    \Delta s(P_j) := \Delta s(u^L, u^M, u^R),
\end{equation*}
where
\begin{equation}
\label{eq:change-strength-scalar}
     \Delta s(u^L, u^M, u^R):=  |u^L - u^M| + |u^M - u^R| - |u^L - u^R| 
\end{equation}
i.e. the (local) difference of total variation before and after the interaction. Notice that, in the scalar case,  in same sign interactions,  $\Delta s(P_j) = 0$, while, in opposite sign interactions, $\Delta s(P_j) \geqslant 0$  (this is  not anymore true in the case of systems).

The \emph{change of speed} across the interaction is the quantity
\begin{equation*}
    \Delta \sigma(P_j) := \Delta \sigma(u^L, u^M, u^R),
\end{equation*}
where $\Delta \sigma(u^L, u^M, u^R)$ is defined differently in the case of same sign interaction and opposite sign interaction, as follows: 
\begin{equation}
\label{eq:change-speed-scalar}
    \Delta \sigma(u^L, u^M, u^R):=
        \norm {\sigma^{\text{after}}-\sigma^{\text{before}}}_{L^1([u^L, u^R])},\\
\end{equation}
where
\begin{align*}
    &\sigma^{\text{after}} =
        D\, {\conv_{\mathclap{[u^L, u^R]}}f}, \quad \sigma^{\text{before}} =
        \begin{cases}
        D \conv\limits_{\mathclap{[u^L, u^M]}}f\, \cup\, D \conv\limits_{\mathclap{[u^M, u^R]}}f, & u^L\leqslant u^M \leqslant u^R,\\
        D \conv\limits_{\mathclap{[u^L,u^M]}} f\,\Big|_{[u^L, u^R]}, & u^L \leqslant u^R \leqslant u^M. 
        \end{cases}
\end{align*}

Notice that, when solving a Riemann problem (e.g. between $u^L$ and $u^R$), the derivative $D \conv_{[u^L, u^R]}f$ of the convex envelope of $f$ (or concave, in case of negative jumps) corresponds to the speeds of the various jumps, rarefactions and contact discontinuities the solution consists of. Hence, $\Delta \sigma(u^L, u^M, u^R)$ can be though as the $L^1$ distance  between the ``speed'' functions. 

\begin{example}
  Assume that at point $P_j$ two (positive) shocks interact, the first one connecting $u^L$ to $u^M$ and traveling with speed $\sigma'$, the second one connecting $u^M$ to $u^R$
  and traveling with speed $\sigma''$. At the interaction point the two shocks merge and produce a new shock connecting $u^L$ to $u^R$ and traveling with speed $\sigma$. This new shock can be thought of as the superposition of the two incoming shocks $(u^L, u^M)$ and $(u^M, u^R)$, which change their speed from $\sigma'$ (resp. $\sigma''$) to $\sigma$. The change of speed is given by
\begin{equation*}
    \Delta \sigma(P_j) := \Delta \sigma(u^L, u^M, u^R) = |\sigma' - \sigma|\cdot |u^M - u^L| + |\sigma'' - \sigma|\cdot |u^R - u^L|. 
\end{equation*}
\end{example}

\subsection{Glimm's program}
\label{ss:glimm-program}

In his seminal paper~\cite{Glimm1965}, Glimm developed a general strategy to prove the estimates in~\eqref{eq:intro-interaction-estimate}. Such strategy may be applied to every scheme (e.g. Glimm scheme or wavefront tracking) where approximate solutions are constructed by patching together local solutions to Riemann problems and taking care of interactions. We briefly describe this strategy.

An \emph{interaction} can be defined in  very general terms (independent of the particular approximation scheme), as follows: there are two maps $u = u(x), \hat u = \hat u(x) \in \mcD$ (we can think to $u$ and $\hat u$ as the approximate solution respectively \emph{before} and \emph{after} the interaction). The two maps $u, \hat u$ coincide outside an interval $(\alpha,\beta)$. On  $(\alpha,\beta)$, both $u, \hat u$ are piecewise constant. The map $u$  has exactly two jumps on $(\alpha, \beta)$, namely between $u^L$ and $u^M$ and between $u^M$ and $u^R$, while the map $\hat u$ has a single jump between $u^L, u^R$, see Fig.~\ref{fig:thm:interaction-estimate} on page~\pageref{fig:thm:interaction-estimate}.

As mentioned before, to each \emph{interaction}, we can associate the change of strength $\Delta s(u^L, u^M, u^R)$ and change of speed $\Delta \sigma (u^L, u^M, u^R)$ of wavefronts involved in the interaction (as we did above~-- in the scalar case~-- in~\eqref{eq:change-strength-scalar},~\eqref{eq:change-speed-scalar}). 

Glimm's strategy to prove~\eqref{eq:intro-interaction-estimate} consists in finding an \emph{interaction functional} \linebreak $Q\colon\mcD \to \R $ with the following properties:
\begin{enumerate}[(i)]
    \item \label{pt:bounds} For each $v \in \mcD$, $Q(v)$ is nonnegative and it is bounded above by a constant depending on $f$ and on $\TV(v)$:
    \begin{equation*}
        0 \leqslant Q(v) \leqslant C(f, \TV(v)).
    \end{equation*}
    \item \label{pt:int-est} In each interaction (using the same notation as above), $Q$ decreases and its decrease bounds the change of strength
    \begin{subequations}
    \begin{equation}
        \Delta s(u^L, u^M, u^R) \leqslant Q(u) - Q(\hat u)  
    \end{equation}
    or the change of speed
    \begin{equation}
    \label{eq:change-speed}
        \Delta \sigma(u^L, u^M, u^R) \leqslant Q(u) - Q(\hat u), 
    \end{equation}
    \end{subequations}
    of the two wavefronts $(u^L, u^M)$, $(u^M, u^R)$ involved in the interaction, depending on which of the two estimates~\eqref{eq:intro-strength}--\eqref{eq:intro-speed} one wants to prove.
    \item \label{pt:no-interaction} In an approximate solutions $u^\e$, the map $t \mapsto Q(u^\e(t))$ is constant as long as no interaction takes place.\footnote{
    Point~\eqref{pt:no-interaction} is actually a consequence of $Q$ decrease in interactions: indeed, shifting a jump of piecewise constant function can be considered as an interaction of this jump with zero jump.
    } 
\end{enumerate}

Notice that the three items above (with~\eqref{eq:change-speed} in item~\eqref{pt:int-est}) correspond to Points~\eqref{thm-main:p1}--\eqref{thm-main:p3} in the~\ref{thm-main} on page \pageref{thm-main}.

Once such an interaction functional is found, it is standard to see, by a telescopic sum argument, that, given an approximate solution $u^\e(t,x)$ (constructed e.g. by wavefront tracking or the Glimm scheme), the map $t \mapsto Q(u^\e(t))$ is non-increasing and the total change of strengths $\sum_j \Delta s(P_j)$ (or the total change of speeds $\sum_j \Delta \sigma(P_j)$) is bounded by
$Q(0)= Q(u^\e(0))\leqslant  C(f, \TV(u^\e(0)) \leqslant C(f,\allowbreak \TV(\bar u))$, hence proving~\eqref{eq:intro-interaction-estimate}.

Glimm's strategy was successful in proving the estimates on the change of \linebreak strength \eqref{eq:intro-strength} and on the change of speed~\eqref{eq:intro-speed} in the GNL/LD case (\cite{Glimm1965,Liu1977}, see also the book \cite{Bressan2000} and references therein) and to prove the estimate on the change of strength~\eqref{eq:intro-strength} in the general case (where no GNL/LD assumption is made), see~\cite{Bianchini2002,Bianchini2003} and also~\cite{BianchiniBressan2005}, where the estimate on the change of strength was proven for vanishing viscosity approximations. 

The estimate for the change of speed~\eqref{eq:intro-speed} in the general case (no GNL/LD assumption) is somewhat a  different story.

\subsection{The change of speed in the general case}\label{ss:CSGC}

The estimate for the change of speeds~\eqref{eq:intro-speed}
was proven in~\cite{BM-CMP, M-Forward} (see also \cite{BM-Glimm, BM-CaseStudy, BM-scalar}) for approximate solutions to~\eqref{eq:syst-cons-law} constructed by the Glimm scheme (but the same proof would work for wavefront tracking solutions), with a strategy similar to Glimm's one, described in the previous sections. 

However, in the mentioned works~\cite{BM-CMP, M-Forward}, the authors do not construct an interaction functional $Q\colon \mcD \to \R$ with the properties~\eqref{pt:bounds}--\eqref{pt:no-interaction} above. Rather, for any \emph{given} approximate solution $u^\e(t,x)$, the authors construct a functional $\mathfrak Q = \mathfrak Q(t)$, with properties analog to~\eqref{pt:bounds}--\eqref{pt:no-interaction}  above: in particular, $\mathfrak Q(0) \leqslant C(f, \TV(\bar u))$, at each interaction time $\Delta \sigma(P_j) \leqslant \mathfrak Q(t_j^-) - \mathfrak Q(t_j^+)$, and $\mathfrak Q$ is constant when no interaction takes place, so that~\eqref{eq:intro-speed} follows. 

The difference w.r.t. Glimm's strategy described in Section~\ref{ss:glimm-program} is that such functional $\mathfrak Q(t)$ of~\cite{BM-CMP, M-Forward}  can not be written in the form $Q(u^\e(t))$ for some $Q\colon \mcD \to \R$, or, in other words, the definition of $\mathfrak Q(t)$ is highly \emph{nonlocal in time}. The reason is that the value of the functional $\mathfrak Q$ at time $t$ depends not only on the map $u^\e(t, \cdot) \in \mcD$, but on the whole approximate solution $u^\e$ on $[0,+\infty) \times \R$ through a nontrivial wave tracing algorithm, which takes care, at any time, of whole past and future history of waves.\footnote{To be precise, in \cite{BM-CMP}, $\mathfrak Q(t)$ depends on $u^\e|_{[0,t] \times \R}$, while in \cite{M-Forward} it depends on $u^\e|_{[t, + \infty) \times \R}$.} In particular, $\mathfrak Q(t)$ depends on the approximation parameter $\e$ (though it enjoys a uniform bound, independent of $\e$)  and on the order in which interaction takes place in the approximate solution $u^\e$. 

This is clearly not an issue, if the goal is just to prove~\eqref{eq:intro-speed}, as it was in~\cite{BM-CMP, M-Forward}. However, the nonlocality in time of $\mathfrak Q$ has also several disadvantages, that we will list and discuss in the next Section~\ref{ss:relevance}.

\medskip

The goal of this paper is precisely to fill this gap and to construct a functional 
\begin{equation*}
    \Upsilon\colon\mcD \to \R
\end{equation*}
with Properties~\eqref{pt:bounds},~\eqref{pt:int-est} (with the estimate~\eqref{eq:change-speed} for the change of speeds) and~\eqref{pt:no-interaction} above, which is completely local in time, meaning, in particular, that $\Upsilon(u)$ is defined for every $u \in \mcD$ (not only for piecewise constant maps) and  its definition is not linked to a specific approximation scheme.

The fact that $\Upsilon$ is purely local in time, in analogy to the classical functionals of Glimm and of Bianchini--Bressan (\cite{Glimm1965, Bianchini2002, Bianchini2003}) and differently than the functionals used to prove~\eqref{eq:change-speed} in~\cite{BM-CMP, M-Forward}, is the most important contribution of this paper.

\subsection{Why is a \emph{local-in-time} functional relevant?}
\label{ss:relevance}

The \emph{locality in time} of $\Upsilon$ has several advantages.

First of all, the definition of $\Upsilon$ in~\eqref{UpsilonScalar} is much cleaner and simpler than the one of the functionals in~\cite{BM-CMP, M-Forward}, where (even in the scalar case~\cite{BM-scalar, M-Forward}) a preliminary construction of a detailed wave tracing  algorithm (also called \emph{enumeration of waves} or \emph{Lagrangian representation} in the mentioned works) was needed. 

Second, our new functional $\Upsilon$ is defined on every  map in the semigroup domain~$\mathcal{D}$, whereas the functionals in~\cite{BM-CMP, M-Forward} are defined on approximate solutions and thus, in particular, just for piecewise constant maps fitting the setting of Glimm scheme or wavefront tracking. The functionals in~\cite{BM-CMP, M-Forward} do not have a general $\BV$ counterpart, because their definition was based on the wave tracing, and such wave tracing can be explicitly constructed only on approximate solutions, where interactions are concentrated on a finite number of points in the $(t,x)$-plane. Our new functional $\Upsilon$, on the contrary, can be computed along trajectories of the semigroup, and, moreover, it is decreasing along such trajectories. We believe this fact can be further exploit to analyze finer properties of the exact solution to~\eqref{eq:syst-cons-law}--\eqref{eq:initial-datum}.  

From a different point of view, the idea of developing a \emph{Lagrangian representation} for conservation laws (as an extension of the method of characteristic) was pioneered by Dafermos \cite{Dafermos1989} for GNL/LD systems. In more recent years, for scalar conservation laws (also in several space dimensions), Lagrangian representations were developed by Bianchini and Marconi \cite{BiaMar16, BiaMar17}, Bianchini, Bonicatto and Marconi
\cite{BiaBonMar17} and Marconi \cite{Marconi2021, Marconi2022}.

On the other side, for systems of conservation laws, a preliminary version of a Lagrangian representation was proposed in~\cite{M-Thesis, M-HYP2016}, but the analysis there was quite involved also because of the absence of a local-in-time functional bounding the change of speeds in interactions. Our new functional $\Upsilon$ should therefore fill this gap and thus it should make easier to construct Lagrangian representation for systems of conservation laws. We plan to develop this program in future works.

A final remark is devoted to viscous systems of conservation laws 
\begin{equation}
    \label{eq:viscous}
    \partial_t u + \partial_x f(u) = \nu \partial_{xx} u, \qquad \nu>0.
\end{equation}
 The Bianchini--Bressan functional \cite{Bianchini2002, Bianchini2003}, defined for the inviscid case, can be easily extened to bound the \emph{change of strength} of ``viscous wavefronts'' in~\eqref{eq:viscous}, as in \cite{BianchiniBressan2005}.

 The functionals in~\cite{BM-CMP, M-Forward}, used to bound the change of speeds in the inviscid case and whose definition was based on a wave tracing algorithm, had clearly very little chances to be extended to viscous conservation laws, since in a solution to~\eqref{eq:viscous} it is not clear how to ``trace'' viscous wavefronts. 

Our new functional $\Upsilon$ does not have an immediate viscous counterpart. However, we believe that the fact that $\Upsilon$ is local-in-time contributes to make clearer why the change of speeds of wavefronts in the inviscid system  stays bounded, and therefore it can provide some insights if one wants to bound the change of speeds of viscous wavefronts (or some related quantities) in a solution to~\eqref{eq:viscous}. Such an estimate could be relevant (but this is nothing more than wishful thinking) in order to prove convergence of vanishing viscosity approximation, when a \emph{physical} viscosity is considered, instead of the \emph{artificial} viscosity $\nu \partial_{xx} u$ considered in~\cite{BianchiniBressan2005}.

\subsection{Structure of the paper and notations}

In Section~\ref{s:def-rig-stat} we give the definition of the new interaction functional $\Upsilon$ (based on the fundamental notion of $\cost$) and we present the rigorous statements corresponding to items~\eqref{thm-main:p1}--\eqref{thm-main:p4} of~\ref{thm-main}.   Section~\ref{s:cost} is the core of the paper, where we prove the fundamental properties of $\cost$. In Section~\ref{s:scalar} we prove items~\eqref{thm-main:p1}--\eqref{thm-main:p2} of~\ref{thm-main}. In Section~\ref{s:const-no-int} we prove item~\eqref{thm-main:p3} of~\ref{thm-main}. Appendix~\ref{s:semi-explicit} contains a formula (Theorem~\ref{th:cost-stress-new}) that expresses $\cost$ in terms of a deterministic procedure. The material of the Appendix is not needed in the rest of the paper, but it provides more insight into the structure of $\cost$.

We use the following notation:
\begin{itemize}
    \item $\conv_{[a,c]} f$ for the convex envelope of the function $f$ over the segment $[a,c]$;
    \item for $f\colon[a,b] \to \R$, $g\colon [b,c] \to \R$, we set 
    \begin{align*}
    \label{eq:cup-def}
        f\cup g\colon [a,c] \to \R,\quad (f\cup g)(y) =
        \begin{cases}
          f(y), &y\in [a,b),\\
          g(y), &y\in [b,c];
        \end{cases} 
    \end{align*}
    \item $u_x$ for the derivative (treated as a measure) of a $\BV$ function $u$;
    \item $\ind A$ for $A\subset\R$;
    \item $\mu\restr{A}\colon B\mapsto \mu(A\cap B)$ the restriction of Radon measure $\mu$ onto a Borel set $A$;
    \item $\sqcup$ for disjoint union;
    \item $g^{\pm }:=\frac12(\abs{g}\pm g)$ for the positive (resp. negative) part of function $g$;
    \item $\partial^\pm$ for right (resp. left) derivative;
    \item $C^{1,1}([a,c])$ for the set of continuously differentiable maps on $[a,c]$ with Lipschitz continuous derivative;
    \item $\osc_{[a,c]}g = \sup_{[a,c]} g - \inf_{[a,c]} g$;
    \item $\psi_\sharp\mu$ for the pushforward of the measure $\mu$ through the map $\psi$.
\end{itemize} 

\section{Definition of \texorpdfstring{$\Upsilon$}{Upsilon} and rigorous statements}
\label{s:def-rig-stat}

From now on, we focus just on the scalar case $n=1$. The analysis of systems ($n>1$) is performed in the companion paper~\cite{BaitMod2}. In particular, we assume that the domain of the semigroup $S_t$ is 
\begin{equation*}
    \mathcal{D} = \BV(\R; \overline \Omega),
\end{equation*}
where $\Omega \subseteq \R$ is open. We assume w.l.o.g. that maps in $\mcD$ are right continuous. 
$\mathcal{D}$ will be as well the domain of our  functional $\Upsilon$. 
We present in this section the precise definition of $\Upsilon$ and we rigorously state its main properties.

\subsection{Definition of  \texorpdfstring{$\cost$}{cost}}
\label{ss:cost-def}

We start with the most important definition of this paper, namely the notion of \emph{cost function}. 

We mostly focus on $C^{1,1}$-regularity of the functions involved, motivated by applications to conservation laws. However, Definitions~\ref{def:H},~\ref{def:mergecost} and~\ref{def:cost} make sense for continuous functions. In Appendix~\ref{s:semi-explicit}, we explore the further structure of $\cost$ in the $C^0$-regularity setting.

\begin{figure}[h]
			\centering
			\includegraphics{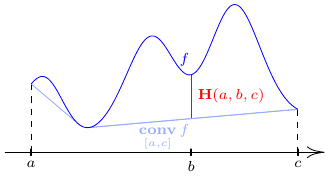}
		\caption{To Def.~\ref{def:H}.}
		\label{fig:Hdef}
	\end{figure}

We need first the definitions of two auxiliary quantities.     
\begin{definition}\label{def:H}
    Let $a \leqslant b\leqslant  c$ and $f\in C^{1,1}([a,c])$. Define the \emph{height} by
    \begin{align}
        \Height_f(a,b,c) =
        f(b) -\Big(\conv_{[a,c]} f\Big)(b).\label{eq:Hdef}
    \end{align}
\end{definition}
    Def.~\ref{def:H} is illustrated in Fig.~\ref{fig:Hdef}.
\begin{remark}
    The quantity $\Height_f(a,b,c)$ is related to (and motivated by) the definition of $\Delta\sigma$ in the scalar case; see~\eqref{eq:def-delta-sigma-scalar-1}--\eqref{eq:def-delta-sigma-scalar-2} below.
\end{remark}

\begin{definition}\label{def:mergecost}
     Let $a\leqslant c$ and $f\in C^{1,1}([a,c])$. Suppose that $b_1,\ldots,b_m \in [a,c]$. We define the \emph{merge cost} of the interval $[a,c]$ w.r.t. to the \emph{mergers} $(b_i)_{i=1}^m$   recursively by
     \begin{align}
         \cost_f([a,c],(b_i)_{i=1}^m) = \begin{cases}
             \Height_f(a,b_1,c),& m = 1,\\
             \cost_f([a,c],(b_i)_{i=1}^{m-1}) +\Height_f(\bleft_m, b_m, \bright_m),& m>1,
         \end{cases}\label{def:mergecost:eq1}
     \end{align}
      where 
     \begin{equation}
    \label{def:mergecost:eq2}
     \begin{aligned}
        \bleft_j& = \max\Big(\{a\}\cup\left( \{b_i\}_{i=1}^{j-1}\cap[a,b_j]\right)\Big),\\  
        \bright_j & = \min\Big(\{c\}\cup\left( \{b_i\}_{i=1}^{j-1}\cap[b_j,c]\right)\Big). 
     \end{aligned}
     \end{equation}
\end{definition}

\begin{figure}[h]
			\centering
			\includegraphics{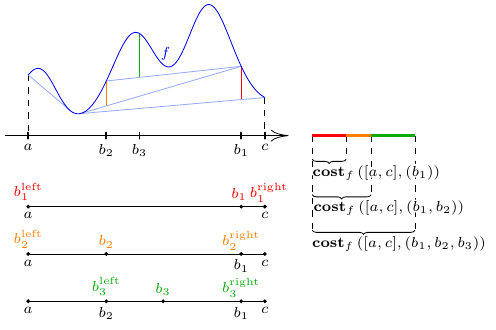}
		\caption{The definition of merge cost.}
		\label{fig:mergercostdef}
\end{figure}

\noindent    Notice that, in general, $\cost_f([a,c],(b_i)_{i=1}^m)$ depends on the order of the points $b_i$. An example is shown in Fig.~\ref{fig:mergercostdef}.
Unwinding the recursion, we also obtain a more explicit expression
\begin{align}
\label{eq:mcostdef}
     \cost_f([a,c],(b_i)_{i=1}^m) = \sum_{i=1}^m\Height_f(\bleft_i, b_i, \bright_i).
\end{align}

We can finally introduce the \emph{cost function} $\cost$.

\begin{definition}\label{def:cost}
        Let $a\leqslant c$ and $f\in C^{1,1}([a,c])$. We define the \emph{cost function} as 
        \begin{align}
        \label{def:costfac:eq}
            \cost_f(a,c) =\sup\Big\{ \cost_f([a,c],(b_i)_{i=1}^m)  \colon \,  m\in \BN, \, (b_i)_{i=1}^m \subset [a,c] \Big\}.
        \end{align}
        Define as well $
            \cost_f(c,a) = \cost_{-f}(a,c)$.
    \end{definition}

\begin{figure}[h]
		\centering
		\begin{minipage}{0.45\textwidth}
			\centering
			\includegraphics{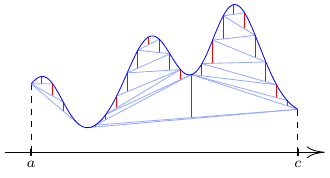}
            \\
			(a) 
		\end{minipage}
        \hfill
        \begin{minipage}{0.45\textwidth}
			\centering
			\includegraphics{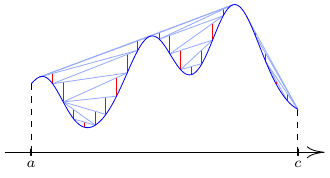}
            \\
			(b) 
		\end{minipage}
		\caption{Quasi-optimal mergers: (a) for $\cost_f(a,c)$; (b) for $\cost_f(c,a)$.}
		\label{fig:quasiopt_mergers}
	 \end{figure}

\noindent We will prove (Prop.~\ref{prop:cost_Linfty_est})  that 
\begin{equation}
    \label{eq:quadr-bd}
    \cost_f(a,c) \leqslant \frac{1}{4} \left\| \frac{\partial^2 f}{\partial \tau^2} \right\|_{L^\infty} |c-a|^2,
\end{equation}
and thus, in particular, $\cost_f(a,b)$ is always a  finite quantity. For the sake of better understanding, we provide an illustration of quasi-optimal configurations  in  Fig.~\ref{fig:quasiopt_mergers}.
\begin{remark}
    In view of conservation laws, Def.~\ref{def:cost} may be treated as counting the total potential change of speed due to same-sign interactions, ``hidden'' inside the jump $[a,c]$.
\end{remark}

\subsection{Definition of the interaction functional \texorpdfstring{$\Upsilon$}{Upsilon}}
\label{ss:intro-intfcn}

We can now define our new interaction functional $\Upsilon$.

\begin{definition}[New interaction functional]
We define $\Upsilon\colon\mcD \to \R_{\geqslant 0}$ by
\begin{align}
\label{UpsilonScalar}
    \Upsilon(u) = \frac14\norm{\derder f}_{L^\infty} \TV(u)^2- \underbrace{\sum_{\substack{ |u_x|(\{x\})>0}} \cost_f(u_-(x), u_{+}(x))}_{=:\Qcost(u)}.
\end{align}
\end{definition}

The fact that $\Upsilon$ is well defined (i.e. it maps $\mcD$ into $\R_{\geqslant 0}$) and that it is a \emph{quadratic} functional (w.r.t. $\TV(u)$) will be shown in Theorem~\ref{thm:quadr-bd}. 

We stress that $\cost_f(u_-(x), u_+(x))$ is a  \emph{purely local} quantity, both in space and in time, meaning that it depends only on the left and the right limits of $u$ at the point $x$,  and on nothing else: in particular, it is \emph{local in space}, as it does not depend on the other jumps or on the continuity points of $u$, and it is \emph{local in time}, as it does not depend on the (approximate) solution starting from $u$.

Our functional $\Upsilon$ aims at bounding the change of speed. The fact that $\Qcost(u)$, i.e. ``cost'' part of $\Upsilon$, is zero if there are no jumps in $u$ is consistent with the fact that, along a solution to the scalar conservation law~\eqref{eq:syst-cons-law}, wavefronts can change speeds only if and where jumps are present. Otherwise, characteristics move with constant speed.

In the next sections we state the main properties of $\Upsilon$, corresponding to properties~\eqref{thm-main:p1}--\eqref{thm-main:p3}, as well as~\eqref{thm-main:p4}--\eqref{thm-main:p4}, listed in the~\ref{thm-main} on page \pageref{thm-main}.

\begin{remark}\label{rmk:quad-cubic}
    The definitions of $\cost_f(a,c)$ and $\Qcost(u)$ have some analogy with the definition of $\mathcal{A}(a,c)$ and $\mathcal{D}(u)$ respectively in \cite[Definitions 2.6, 2.7]{Caravenna2008}. In the present paper, $\cost$ is used to control the variation of speed $\Delta \sigma(P_j)$ (notation as in Section~\ref{ss:basic-pbl}), while $\mathcal{A}(a,c)$ in \cite{Caravenna2008} is used to control the variation of the entropy dissipation $\Delta \mathcal{D}(P_j)$ across interactions.

    The crucial difference is that in \cite{Caravenna2008} the quantity to keep under control, i.e. $\Delta 
\mathcal{D}(P_j)$, is, \emph{by definition}, the difference (before and after the interaction) of a functional associated to $u$ (namely, the total entropy dissipation $\mathcal{D}(u)$).

    In our case, on the contrary,  there is no natural functional such that $\Delta \sigma (P_j)$ can be written as difference of such functional (before and after the interaction).

    The main task accomplished in this paper is precisely to \emph{invent} a functional (namely, $\Qcost$) whose positive variation across same sign interaction provides an upper bound for $\Delta \sigma(P_j)$, whilst its negative variation across opposite sign interactions is bounded by the decrease of $\TV(u)^2$. 
     
\end{remark}

\subsection{Nonnegativity and quadratic bound}\label{ss:non-neg-quadr-bound}

The first result claims that $\Upsilon$ is \emph{quadratic} (w.r.t. $\TV(u)$), and  that it is nonnegative. This corresponds to Point~\eqref{thm-main:p1} in the~\ref{thm-main}. While the upper bound is straightforward, nonnegativity is not obvious, because $\Upsilon$ is the \emph{difference} of two nonnegative quantities.

\begin{theorem}[Nonnegativity and quadratic bound]
\label{thm:quadr-bd}
    For all $u \in \mcD$,
    \begin{equation*}
        0 \leqslant \Upsilon(u) \leqslant\norm{\derder f}_{L^\infty} \frac{1}{4} \TV(u)^2.
    \end{equation*}
\end{theorem}

\noindent The proof is given in Section~\ref{ss:scalar-quadratic-bound}.

\subsection{Interaction estimate}

The second result is the main property of $\Upsilon$, namely the interaction estimate.   It says that the change of wavefront speed across an interaction is bounded by the decrease of $\Upsilon$. It  corresponds  to Point~\eqref{thm-main:p2} in the~\ref{thm-main}, up to a harmless multiplicative factor $\frac{1}{2}$ (see also~\eqref{eq:change-speed}).

\begin{figure}[h]
    \centering
    \includegraphics{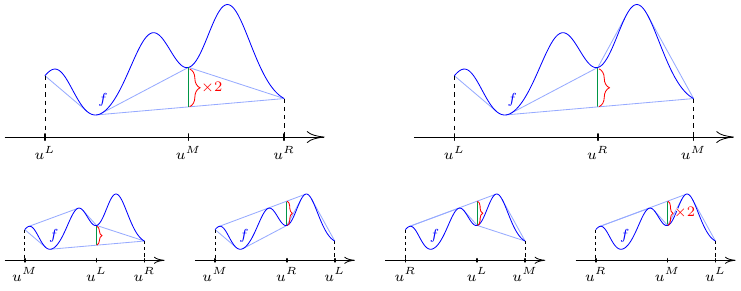}
    \caption{To the definition of $\Delta \sigma$.}
    \label{fig:change_of_speed}
\end{figure}

First of all, we rigorously introduce the quantity $\Delta \sigma(u^L, u^M, u^R)$, which represents the \emph{change of speed} of wavefronts involved in the interaction between the Riemann problems $(u^L, u^M)$ and $(u^M, u^R)$. In order to define it, we have to distinguish six cases, according to the order of the three values $u^L, u^M, u^R$ (see Fig.~\ref{fig:change_of_speed}).
We have already defined $\Delta \sigma(u^L, u^M, u^R)$ in~\eqref{eq:change-speed-scalar} in the  two cases $u^L \leqslant u^M \leqslant u^R$ and $u^L \leqslant u^R \leqslant u^M$. Notice that~\eqref{eq:change-speed-scalar} is equivalent to
\begin{equation}\label{eq:def-delta-sigma-scalar-1}
    \Delta \sigma(u^L, u^M, u^R) = 
    \begin{cases}
        2 \Height_f(u^L, u^M, u^R), & \text{if } u^L \leqslant u^M \leqslant u^R, \\
        \Height_f(u^L, u^R, u^M), & \text{if } u^L \leqslant u^R \leqslant u^M. 
    \end{cases}
\end{equation}
Similarly, we define $\Delta \sigma(u^L, u^M, u^R)$ in the remaining four cases, as follows:
\begin{equation}\label{eq:def-delta-sigma-scalar-2}
    \Delta \sigma(u^L, u^M, u^R) = 
    \begin{cases}
        \Height_f(u^M, u^L, u^R), & \text{if } u^M \leqslant u^L \leqslant u^R, \\
        \Height_{-f}(u^M, u^R, u^L), & \text{if } u^M \leqslant u^R \leqslant u^L, \\ 
        \Height_{-f}(u^R, u^L, u^M), & \text{if } u^R \leqslant u^L \leqslant u^M, \\ 
        2 \Height_{-f}(u^R, u^M, u^L), & \text{if } u^R \leqslant u^M \leqslant u^L. 
    \end{cases}
\end{equation}
Observe that in \emph{same sign interactions} (i.e. when the two jumps $(u^L, u^M)$ and $(u^M, u^R)$ have the same sign, $\Delta \sigma$ is given by twice the height, whereas in \emph{opposite sign interactions} $\Delta \sigma$ is given by the height. This is due, roughly speaking, to the fact that in opposite-sign interaction some ``waves'' are canceled, and thus there are less waves for which one wants to control the change of speed. 

\begin{theorem}[Interaction estimate]
    \label{thm:interaction-estimate}
    Let $u, \hat u \in \mcD$.  Let $I \subseteq \R$ be an open interval. Assume that
    \begin{enumerate}[(a)]
        \item $u(x) = \hat u(x)$ for $x \notin I$,
        \item on $I$ both $u, \hat u$ are piecewise constant,
        \item on $I$, $u$  has exactly two jumps from $u^L$ to $u^M$ and from $u^M$ to $u^R$ respectively, whilst $\hat u$  has exactly one jump from $u^L$ to $u^R$.
    \end{enumerate}
    Then 
    \begin{equation}
    \label{eq:int-est-main}
        \frac{1}{2} \Delta \sigma (u^L, u^M, u^R) \leqslant \Upsilon (u) - \Upsilon(\hat u).
    \end{equation}
    \begin{figure}[h]
        \centering
        \includegraphics{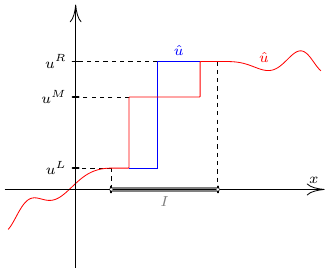}
        \caption{To Theorem~\ref{thm:interaction-estimate}.}
        \label{fig:thm:interaction-estimate}
    \end{figure}
\end{theorem}

\noindent The proof of Theorem~\ref{thm:interaction-estimate} is given in Sections~\ref{ss:scalar-same-sign}--\ref{ss:scalar-opposite-sign}.

\subsection{\texorpdfstring{$\Upsilon$}{Upsilon} is constant if no interaction takes place}

The third property of $\Upsilon$ corresponds to Point~\eqref{thm-main:p3} in the~\ref{thm-main}.
\begin{theorem}[$\Upsilon$ is constant if no interactions]
\label{thm:const-if-no-interaction}
  Let $u \in \mcD$ be piecewise constant with distance between adjacent jump points bounded from below by $L_0>0$. Then   
  for all $t \in [0,t_0)$, $$\Upsilon(S_t u) = \Upsilon(u),$$
  with 
  \begin{equation*}
      t_0 := \frac{L_0}{\osc\limits_{[\inf u, \sup u]} \left( \der f \right)}
  \end{equation*}
  
\end{theorem}
\noindent The proof is given in Section~\ref{s:const-no-int}.

\subsection{Additional properties of \texorpdfstring{$\Upsilon$}{Upsilon}}
\label{ss:additional-properties}

We now list the additional properties of~$\Upsilon$ corresponding to \ref{thm-main}, Point~\eqref{thm-main:p4}. We do not prove these properties in the present paper. The proof (for the more general system case) is given in the companion paper~\cite{BaitMod2}. Reason for this choice is that these are not the main properties of $\Upsilon$, and therefore we find it better, in this case, not to present two separate proofs (for the scalar and the system case). 

We first introduce two operators, $\Cr{}$ and $\RcOp$ on $\BV(\R;\R)$, whose aim is to approximate a given $u \in \BV(\R; \R)$ with piecewise constant functions. $\Cr{}$ removes all total variation outside a given interval. $\RcOp$ collapses the variation of $u$ on given, pairwise disjoint, adjacent intervals $\{I_i\}_{i=1}^m$ into single jumps located at points $x_i \in I_i$.

\begin{definition}
    \label{def:crop}
Let $I \subseteq \R$ be a nonempty interval. 
For $u \in \BV(\R;\R)$, we define the \emph{cropping} $\Cr{I}(u)$ as the map $v \in \BV(\R;\R)$ such that
\begin{align*}
    v_x = v_x\restr{I} =  u_x\restr{I} \quad\text{ and }v =  u\text{ a.e. on } I.
\end{align*}
\end{definition}

Let $I \subseteq \R$ be an interval. We denote by $\mcP(I)$ the set of all finite partitions of $I$ into subintervals.

\begin{definition}
    \label{def:riemann-cut}
Let  $\mcF= \{ I_i\}_{i=1}^m \in \mcP(I)$. We say that the set $X = \{x_i\}_{i=1}^m$ is \emph{compatible} with $\mcF$  if $x_i \in I_i$. 
In that case, for $u \in \BV(\R; \R)$, we define the \emph{Riemann cut} $\Rc\mcF X(u)$ as the map $\hat u$ such that
\begin{equation*}
    \hat u(- \infty) = u(-\infty) \text{ and } \hat u_x = u_x - u_x\restr {I} + \sum_{i=1}^m u_x(I_i) \delta_{x_i}.
\end{equation*}
If  $\mcF = \{I\}$ and $X = \{\hat x\}$, we say that the Riemann cut is \emph{simple} and write $\Rc I{\hat x}(u)$.
\end{definition}
\begin{figure}[h]
    \centering
    \includegraphics{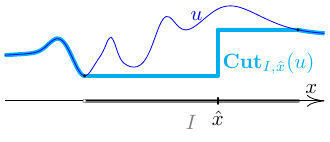}
    \caption{A (simple) Riemann cut.}
    \label{fig:riemann_cut}
\end{figure}
\noindent Notice that, for $u\in\BV(\R;\R)$, the value of $\Upsilon(\Rc \mcF X(u))$ does not depend  on the particular choice of $X$. In view of that, we will write shortly $\Upsilon(\Rc \mcF \bullet(u))$. 

We can now rigorously state the properties corresponding to \ref{thm-main}, Point \eqref{thm-main:p4}.

\begin{theorem}
\label{thm:further-properties}

    The functional $\Upsilon$ has the following additional properties.
    \begin{enumerate}[(a)]
        \item \label{prop:sup-prop-crop} {\underline{Approximation w.r.t. $\Cr{}$}:}
    Let $\Theta \subseteq \R$ be such that $\inf \Theta = - \infty$, $\sup \Theta = + \infty$. For every $u \in \mcD$, 
    \begin{align*}
        \Upsilon(u) = \sup_{\substack{I \subseteq \R \\ \inf I,  \sup I \in \Theta}} \Upsilon(\Cr{I} (u)).
    \end{align*}

    \item \label{thm:sup-prop-cut}
    {\underline{Approximation w.r.t. $\RcOp$}:}
  Let  $\Theta\subseteq I$ be dense in $I$. Then, for every $u\in\mcD$,
\begin{align*}
    \Upsilon(u) = \sup \Big\{ \Upsilon \big(\Rc \mcF \bullet(u) \big) \colon \mcF \in \mcP_\Theta(I) 
    \Big\},
\end{align*}
where 
$\mcP_{\Theta}(I) = \big\{\mcF \in \mcP(I) \text{ s.t. } \forall J \in \mcF, \inf J, \sup J \in \Theta \cup \{ \inf I, \sup I\}\big\}$.
    \item \label{cor:riemann-cut} {\underline{Monotonicity across Riemann cuts}:}
Let  $u\in \mcD$. Let $I\subseteq\R$ be an interval and $\hat x\in I$. Then $\Upsilon(\Rc I{\hat x} (u))\leqslant\Upsilon(u)$.

        \item \label{thm-main:lyap} {\underline{Lyapunov property}:} For every $\bar u \in \mcD$, $t \mapsto \Upsilon(S_t \bar u)$ is non-increasing.

        \item \label{thm-main:lsc} {\underline{Lower semicontinuity}:} If $u_m \to u$ in $\Loneloc$, then $\Upsilon(u) \leqslant \liminf_m \Upsilon(u_m)$.

    \item  \label{thm-main:stab} {\underline{Stability w.r.t. perturbation of state curve}:} Let $u \in \mcD$ and  let $\Phi\colon \overline{\Omega} \to \overline{\Omega}$ be Lipschitz continuous. Then
            \begin{align*}
            \abs{\Upsilon(u)-\Upsilon(\Phi(u))}= 
            \const(1+M^2)M\norm{\Phi-\id}_{W^{1,\infty}},
            \end{align*}
            where
            $M := \max\{\TV(u),\TV(\Phi\circ u)\}$ and $\id(p) \equiv p$.

        \item \label{thm-main:invar} {\underline{Invariance w.r.t. reparametrization of state curve}:} Let $u\in\mcD$ and $\psi\colon \R\to\R$ be strictly increasing. Let $\tilde u\in\mcD$ be such that $\tilde u (-\infty) = u(-\infty)$ and $\tilde u_x = \psi_{\sharp}u_x$. Then $\Upsilon(\tilde u) = \Upsilon(u)$.

    \end{enumerate}
\end{theorem}
\begin{remark}
    We comment on the various properties listed in Theorem~\ref{thm:further-properties}.
    \begin{itemize}
        \item[{\it (a)--(b)}] 
        These fundamental properties claim that the value of $\Upsilon(u)$ can be obtained by approximating $u$ with piecewise constant maps, with variation localized on compact sets. 
        \item[{\it (c)--(e)}] These are common properties of interaction functionals in the literature about hyperbolic conservation laws. We claim that our $\Upsilon$ satisfies them as well.  They are consistent with the fact that ``$\Upsilon$ decreases across interactions'' (Theorem~\ref{thm:interaction-estimate}). 
        \item[{\it (f)--(g)}]  These statements claim that $\Upsilon$ is stable (resp. invariant) under reparameterization of dependent (resp. independent) variable.  
    \end{itemize}
\end{remark}

\section{Cost function}
\label{s:cost}

This section is the core of the paper. We have already defined $\cost$ in Section~\ref{ss:cost-def}. Now we prove its main properties. 
The cost function $\cost$ is the original and most relevant part in the definition of our new interaction functional $\Upsilon$. All properties of $\Upsilon$ stated in the introduction follow from the properties of $\cost$ proven in this section.

In Section~\ref{subsec:bounds}, we prove bounds on $\cost$, which will be the main ingredient of the proof of Theorem~\ref{thm:quadr-bd}. In Section~\ref{ss:cost-interaction}, we prove a superadditive estimate and a subadditive one for $\cost$, which will be the fundamental tool of the proof of the interaction estimate in Theorem~\ref{thm:interaction-estimate}. In Section~\ref{ss:cost-persistency}, we prove \emph{persistency} of $\cost$, which will be the core of the proof of Theorem~\ref{thm:const-if-no-interaction}. In Section~\ref{ss:cost-stability} we prove that the definition of $\cost_f (a,c)$ is stable under perturbation of the function $f$ as well as of  the interval $[a,c]$. This last property is not used in this paper, but it is fundamental in the analysis for systems in the companion paper \cite{BaitMod2}. Finally, in Section~\ref{ss:lsc-sup}, we show an interesting additional property of $\cost$ (actually not needed in the proof of \ref{thm-main}), namely its l.s.c. w.r.t. uniform convergence of the fluxes.

\subsection{Bounds on \texorpdfstring{$\cost$}{cost}}\label{subsec:bounds}

In this section we prove a bound on $\cost_f(a,c)$ in terms of the length of the interval $[a,c]$, and the $L^1$ norm (resp. $L^\infty$ norm) of the second derivative of $f$ in Prop.~\ref{prop:cost_Lone_est} (resp. Prop.~\ref{prop:cost_Linfty_est}). The estimate in terms of the $L^1$ norm is needed for the stability properties in Section~\ref{ss:cost-stability} below, whereas the estimate in terms of $L^\infty$ norm is the core of the proof of Theorem~\ref{thm:quadr-bd}.  

\smallskip

\begin{figure}[h]
    \centering
    \begin{minipage}{0.45\linewidth}
    \centering
    \includegraphics{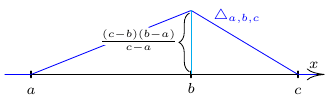}  
    \\
    (a)
    \end{minipage}
        \hfill
    \begin{minipage}{0.45\linewidth}
    \centering
    \includegraphics{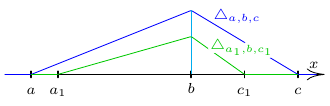}  
    \\
    (b)
    \end{minipage}
    \caption{(a) the function $\triangle_{a,b,c}$; (b)  to the proof of Prop.~\ref{prop:Hbounds}.}
    \label{fig:triangle}
\end{figure}

We start by introducing, for $a\leqslant b\leqslant c$, the auxiliary function
\begin{align}
            \triangle_{a,b,c} (x) = 
			\begin{cases}
				\frac{(c-b)(x-a)}{c-a},& x\in[a,b],\\
				\frac{(b-a)(c-x)}{c-a},& x\in[b,c],\\
                0, & x\notin [a,c],
			\end{cases}\label{eq:triangle}
        \end{align}
see Fig.~\ref{fig:triangle}(a), which we use now to bound $\Height(a,b,c)$. (In case $a = c$, we put $\triangle_{a,b,c} = 0$.)
\begin{prop}\label{prop:Hbounds}
    If $a\leqslant b\leqslant c$ and $f\in C^{1,1}([a,c])$, then
		\begin{align}
		    \Height_f(a,b,c) \leqslant 	\int_a^c \triangle_{a,b,c}(x)\left(\derder[x] f\right)^-(x)\, dx\leqslant\frac12(b-a)(c-b)\norm{\left(\derder[x] f\right)^-}_{L^\infty}.\label{prop:Hbounds:eq1}
		\end{align}
	\end{prop}
	\begin{proof} 
        The second inequality is a straightforward consequence of the identity 
        \begin{align*}
            \norm{\triangle_{a,b,c}}_{L^1} = \frac12(b-a)(c-b).
        \end{align*}
        Let us prove the first one.    
        Assume that $\Height_f(a,b,c)>0$ (i.e. $f(b) > \conv_{[a,c]}f(b)$), otherwise the statement is trivial.
        
		Consider first the case in which $f > \conv_{[a,c]} f$ on the whole interval $(a,c)$. In this case
		\begin{align*}
			&\Height_f(a,b,c) 
             = f(b) - \conv_{[a,c]}f(b) = f(b) - \left[f(a) + \frac{f(c)-f(a)}{c-a}(b-a) \right] \\
            &= \frac{c-b}{c-a}\int_{a}^b 
            \der[x] f
            \, dx - \frac{b-a}{c-a} \int_b^c \der[x] f\, dx \\
			&= \int_{a}^c \der[x]{\triangle_{a,b,c}}{\der[x]f}\, dx
			= - \int_{a}^c{\triangle_{a,b,c}}(x){\derder[x]f}\,dx\leqslant  
            \int_a^c{\triangle_{a,b,c}(x)}{\left(\derder[x] f\right)^-}(x)\, dx,
		\end{align*}
        where in the second last equality we used that the derivative of $\Delta_{a,b,c}$ on $(a,b)$ is constantly equal to $(c-b)/(c-a)$, and similarly it is constantly equal to $-(b-a)/(c-a)$ on $(b,c)$, while in the last equality we integrated by parts. 
        
        Suppose now that $(a_1,c_1)\subset[a,c]$ is a maximal interval containing $b$ such that $f > \conv_{[a,c]}f$ on $(a_1,c_1)$. Then
        \begin{gather*}
            \Height_f(a,b,c) =\Height_f(a_1,b,c_1) \leqslant 	\int_{a}^{c}{\triangle_{a_1,b,c_1}}(x){ \left(\derder[x] f\right)^-}(x)\, dx.
        \end{gather*}
        as we just have shown. Since $a\leqslant a_1\leqslant c_1\leqslant c$, it follows that $\triangle_{a_1,b,c_1}(x)\leqslant\triangle_{a,b,c}(x)$ (see Fig.~\ref{fig:triangle}(b)), and the proof is complete.
	\end{proof}
    Now we are to exploit Prop.~\ref{prop:Hbounds} to develop a bound on $\cost$ functional. The strategy is to mimic the construction of merge cost in Def.~\ref{def:mergecost} using the functions $\triangle_{a,b,c}$ instead of the quantities $\Height_f(a,b,c)$, and derive a pointwise bound for the total sum of them.

\begin{figure}[h]
		\centering
        \begin{minipage}{0.45\textwidth}
			\centering
			\includegraphics{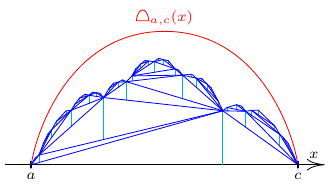}
            \\
			(a) 
		\end{minipage}
		\begin{minipage}{0.45\textwidth}
			\centering
			\includegraphics{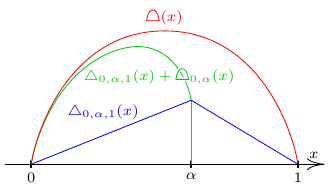}
            \\
			(b) 
		\end{minipage}
        \hfill
        
		\caption{To Lemma~\ref{lemma:psi_triangle}.}
		\label{fig:psi}
    
	\end{figure}

    Denote
    \begin{align*}
        \psishape(x) =\begin{cases}
			x - x\log(2x), &x\in\left(0,\frac12\right],\\ 
			\psishape(1-x), &x\in\left[\frac12, 1\right),\\
            0,& x\notin (0,1),
		\end{cases}
		\qquad 
		\psishape_{a,c}(x) = (c-a)\cdot\psishape\left(\frac{x-a}{c-a}\right).
    \end{align*}
    \begin{lemma}\label{lemma:psi_triangle}
        Let $a,c,(b_i)_{i=1}^m,(\bleft_i)_{i=1}^m, (\bright_i)_{i=1}^m$ be as in Def.~\ref{def:mergecost}. Then for every $x\in[a,c]$,
        \begin{align}
            \sum_{i=1}^m\triangle_{\bleft_i,b_i,\bright_i}(x)\leqslant \psishape_{a,c}(x)\label{lemma:psi_triangle:eq1}
        \end{align}
        (see Fig.~\ref{fig:psi}(a) for an example).
    \end{lemma}
    
\begin{figure}[h]
		\centering
        \begin{minipage}{0.45\textwidth}
			\centering            \includegraphics{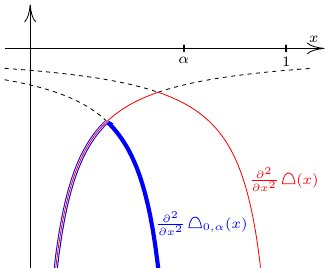}\\
            (a)
            \end{minipage}
            \hfill
            \begin{minipage}{0.45\linewidth}
                \centering 
                \vspace{8mm}\includegraphics{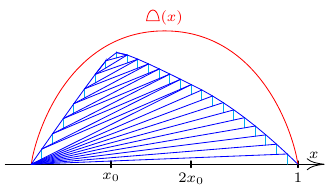}\vspace{6mm}\\
                (b)
            \end{minipage}
		\caption{(a) to Lemma~\ref{lemma:psi_triangle}; (b) to Remark~\ref{remark:psi_triangle_sharpness}.}
		\label{fig:Psi_triangle_ineq_sharpness}
	\end{figure}

    \begin{proof}
        Take $\alpha\in(0,1)$.
        One has
        \begin{align*}
			\derder[x]{\psishape(x)} = \begin{cases}
				-\frac1x,& x\in(0,\frac12],\\
				-\frac1{1-x}, &x\in[\frac12, 1), 
			\end{cases}\qquad
			\derder[x]{\psishape_{0,\alpha}(x)} = \begin{cases}
				-\frac1x,& x\in(0,\frac\alpha2],\\
				-\frac1{\alpha-x}, &x\in[\frac\alpha2, 1),
			\end{cases}
		\end{align*}
        see Fig.~\ref{fig:Psi_triangle_ineq_sharpness}(a),
        so that the function $\psishape - \psishape_{0,\alpha}$ is convex on $[0,\alpha]$. Moreover,  at $x=0$  the value of $\psishape - \psishape_{0,\alpha}$ equals $0$ and its right derivative equals $\log\frac1\alpha>1-\alpha$,
        which implies 
        \begin{align*}
            \psishape_{0,\alpha}(x) + \triangle_{0,\alpha,1}(x)\leqslant\psishape(x), \quad 0\leqslant x\leqslant \alpha,
        \end{align*}
        see Fig.~\ref{fig:psi}(b). By symmetry and scaling, we obtain for every $\alpha<\beta<\gamma$ that 
        \begin{align*}
            \triangle_{\alpha,\beta,\gamma}(x)\leqslant\psishape_{\alpha,\gamma}(x)-\psishape_{\alpha,\beta}(x)-\psishape_{\beta,\gamma}(x),\quad x\in\BR.
        \end{align*}
        
        Applying this inequality to $(\alpha,\beta,\gamma) = (\bleft_i, b_i,\bright_i)$, $i = 1,\ldots, m$, we estimate the l.h.s. of~\eqref{lemma:psi_triangle:eq1} by a telescopic sum that reduces to 
        \begin{gather*}
            \psishape_{a,c}(x) - \psishape_{a,b_{(1)}}(x) - \psishape_{b_{(1)}, b_{(2)}}(x) -\ldots - \psishape_{b_{(m)}, c}(x) \leqslant \psishape_{a,c}(x),
        \end{gather*}
        where $b_{(1)}\leqslant\ldots\leqslant b_{(m)}$ denotes the rearrangement of $b_1,\ldots,b_m$.
    \end{proof}
    \begin{remark}\label{remark:psi_triangle_sharpness}
        For a fixed $x= x_0\in[a,c]$, the estimate~\eqref{lemma:psi_triangle:eq1} is asymptotically sharp. For the case $a=0, c = 1$, $x_0\leqslant\frac12$, the quasi-optimal sequence of subdivisions is illustrated in Fig.~\ref{fig:Psi_triangle_ineq_sharpness}(b). In the limit, this construction naturally leads to the function~$\psishape$.
    \end{remark}
    Now we are in a position to establish the $L^1$ bound for $\cost$.
    \begin{prop}\label{prop:cost_Lone_est}
        Let $a<c$ and $f\in C^{1,1}([a,c])$, then
        \begin{align}
            \cost_f(a,c)\leqslant\int_a^c{\psishape_{a,c}}(x){\left(\derder[x]{f}\right)^-}(x)\, dx\leqslant\frac12(c-a)\norm{\left(\derder[x]{f}\right)^-}_{L^1}.\label{prop:cost_Lone_est:eq1}
        \end{align}
    \end{prop}
    \begin{proof}
        Apply the alternative definition of merge cost in~\eqref{eq:mcostdef} and the first inequality from Prop.~\ref{prop:Hbounds} to get
    \begin{equation*}
        \begin{aligned}
            \cost_f([a,c],(b_i)_{i=1}^m) 
            & = \sum_{i=1}^m\Height_f(\bleft_i, b_i, \bright_i)  \leqslant \sum_{i=1}^m 
            \hspace{4mm}
            \mathclap{\int_{\bleft_i}^{\bright_i}} 
            \hspace{19mm}
            \mathclap{\triangle_{\bleft_i, b_i, \bright_i}(x)\left(\derder[x] f\right)^-}
            \hspace{16.5mm}
            (x)\, dx \\
            & = \int_a^c \sum_{i=1}^m \triangle_{\bleft_i, b_i, \bright_i}(x) \left(\derder[x] f\right)^-(x)\, dx.
        \end{aligned}
    \end{equation*}
    Now apply Lemma~\ref{lemma:psi_triangle} to the integrand,  and then the identity $\norm{\psishape_{a,c}}_{L^\infty} = \frac{c-a}2$.
    \end{proof}
    We also provide a direct estimate in terms of the $L^\infty$-norm of second derivative.
    \begin{prop}\label{prop:cost_Linfty_est}
        Let $a<c$ and $f\in C^{1,1}([a,c])$, then
        \begin{align}
            \cost_f(a,c)\leqslant\frac14(c-a)^2\norm{\left(\derder[x]{f}\right)^-}_{L^\infty}.\label{prop:cost_Linfty_est:eq1}
        \end{align}
    \end{prop}
        \begin{proof}
            Let $a,c,\{b_i\}_{i=1}^m,\{\bleft_i\}_{i=1}^m, \{\bright_i\}_{i=1}^m $ be as in Def.~\ref{def:mergecost}. Then by the second inequality in Prop.~\ref{prop:Hbounds}, we have
            \begin{align*}
                 \cost_f([a,c],(b_i)_{i=1}^m) \leqslant\frac12\norm{\left(\derder[x]{f}\right)^-}_{L^\infty}\sum_{i=1}^m (b_i-b_i^{\rm left})(b_i^{\rm right}-b_i).
            \end{align*}
            Rewriting the $i$-th term of the sum as $\frac12\big[(\bright_i-\bleft_i)^2-(b_i-\bleft_i)^2-(\bright_i-b_i)^2\big]$, we obtain a telescopic sum bounded by $\frac12(c-a)^2$ from above, similar to the proof of Lemma~\ref{lemma:psi_triangle}.
        \end{proof}
        \begin{remark}\label{rmk:prec-consts}
        \begin{enumerate}
            \item    
            One could derive an analogue of Prop.~\ref{prop:cost_Linfty_est} with worse constant ($\frac12$ instead of $\frac14$) from
            Prop.~\ref{prop:cost_Lone_est}. 
            \item The constants in Prop.~\ref{prop:cost_Lone_est} and Prop.~\ref{prop:cost_Linfty_est} are optimal. Indeed, consider
            \begin{align*}
                g(x) = -x^2,\quad h_n(x) = \begin{cases}
                    -\abs{x}, &x\notin\left(-\frac1n,\frac1n\right),\\
                    -\frac{nx^2}{2}-\frac1{2n},&x\in\left(-\frac1n,\frac1n\right),
                \end{cases} 
            \end{align*}
            then, using Theorem~\ref{th:cost-stress-new} in the Appendix, one obtains
            \begin{align*}
                &\cost_g({-1,1}) = 2, \hspace{12.8mm}\norm{{\derder[x]g}}_{L^\infty([-1,1])} \hspace{.8mm}= 2,\\
                & \cost_{h_n}(-1,1) = 2-\frac1n,\quad\norm{{\derder[x]{h_n}}}_{L^1([-1,1])} = 2,
            \end{align*}
            so that Prop.~\ref{prop:cost_Linfty_est} is fulfilled and Prop.~\ref{prop:cost_Lone_est} is asymptotically fulfilled as $n\to\infty$.
        \end{enumerate}
        \end{remark}

\subsection{Interaction-like estimates}
\label{ss:cost-interaction}

We prove  a superadditive and subadditive estimate for $\cost_f(a,c)$. 
We call these bounds \emph{interaction-like estimates} because we will use
\begin{itemize}
    \item the superadditive bound (Prop.~\ref{prop:superadd}) to claim that, in a \emph{same-sign} interaction, the decrease of the functional $\Upsilon$ controls the change of speed of wavefronts involved in the interaction;
    \item  the subadditive bound (Prop.~\ref{prop:subadd}) to claim that, in a \emph{opposite-sign} interaction, $\Qcost$ may increase, but its increase is controlled (up to a multiplicative constant)  by the amount of canceled waves, and thus by the decrease of $\TV^2$, so that overall $\Upsilon$ decreases also across opposite sign interactions. 
\end{itemize}
This section contains therefore already the ideas needed to prove Theorem~\ref{thm:interaction-estimate}.

	\begin{prop}\label{prop:superadd}
		If $a\leqslant b\leqslant c$ and $f\in C^{1,1}([a,c])$, then
		\begin{equation}
        \cost_f(a,c)\geqslant\cost_f(a,b)+\cost_f(b,c)+\Height_f(a,b,c).\label{prop:superadd:eq1}
		\end{equation}
	\end{prop}
    \begin{proof}
        Take any mergers resulting in $[a,b]$ and $[b,c]$, unite them, add first merger $b$, and finally apply the definition of $\cost_f(a,c)$.
    \end{proof}
    Now we prove a sub-additive estimate.
	\begin{prop}\label{prop:subadd}
		If $a\leqslant b\leqslant
        c$ and $f\in C^{1,1}([a,c])$, then
		\begin{equation}
        \cost_f(a,c)\leqslant\cost_f(a,b)+\cost_f(b,c)+
        \frac12
            (b-a)(c-b)\norm{\left(\derder[x]f\right)^-}_{L^\infty}
        .
        \label{prop:subadd:eq1}
		\end{equation}
	\end{prop}
    For the proof of Prop.~\ref{prop:subadd} we need the following  lemma illustrated in Fig.~\ref{fig:lemma:xyzt}.
    
    \begin{lemma}\label{lemma:xyzt}
        If $a\leqslant b_1\leqslant b\leqslant c$, $a<b$ and $f\in C^{1,1}([a,c])$ then
        \begin{align*}
            \Height_f(a,b_1,c)-\Height_f(a,b_1,b)&\leqslant\frac{b_1-a}{b-a}\Height_f(a,b,c) .
        \end{align*}
    \end{lemma}
    \begin{proof}
    \begin{figure}[h]
        \centering
        \includegraphics{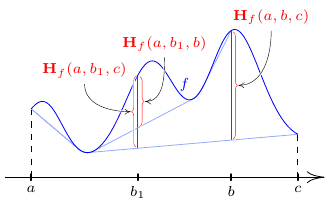}
        \caption{To Lemma~\ref{lemma:xyzt}.}
        \label{fig:lemma:xyzt}
    \end{figure}
    
    It is not difficult to see that the function $g = \conv_{[a,b]}f-\conv_{[a,c]}f$ is convex on $[a,b]$ (for a proof see, for instance,~\cite[Prop.~A.4]{BM-CMP}). 
        Then,
        \begin{align*}
            \Height_f(a,b_1,c)-\Height_f(a,b_1,b) = g(b_1)\leqslant\frac{b-b_1}{b-a}g(a) +\frac{b_1-a}{b-a}g(b).
        \end{align*}
        Now observe that $g(a) = 0$ and $g(b) = \Height_f(a,b,c)$.         
    \end{proof}
    \begin{proof}[Proof of Prop.~\ref{prop:subadd}]

    Let $\{b_i\}_{i=1}^m\subset [a,c]$ be a sequence of points.
    We argue by induction on $m$ that 
    \begin{align}
        &\cost_f([a,c],(b_i)_{i=1}^m)\label{prop:subadd:eq2}\\
        \leqslant&\cost_f(a,b)+\cost_f(b,c)  
         +\frac12
            (b-a)(c-b)\norm{\left(\derder[x]f\right)^-}_{L^\infty}.\nonumber
    \end{align}
    The conclusion clearly follows from~\eqref{prop:subadd:eq2}.

    \textit{The base of induction, $m = 0$}. There is nothing to prove.

    \textit{The inductive step.}
    W.l.o.g., assume that $a<b<c$ and  $b_1\in(a,b]$. Let $(b_i')_{i=1}^{m'}$ and $(b_i'')_{i=1}^{m''}$ be the subsequences of $(b_i)_{i=2}^m$ consisting of points which are less than $b_1$ and not less than $b_1$, respectively. Then,
    \begin{align*}
        &\cost_f([a,c],(b_i)_{i=1}^m) =\Height(a,b_1,c) + \cost_f([a,b_1],(b_i')_{i=1}^{m'})+\cost_f([b_1,c],(b_i'')_{i=1}^{m''}).
    \end{align*}
    Let us bound the first term in r.h.s. by Lemma~\ref{lemma:xyzt}; the second, by $\cost_f(a,b_1)$; the third, by the inductive assumption. Therefore,
    \begin{align*}
        \cost_f([a,c],(b_i)_{i=1}^m) \leqslant&\Height_f(a,b_1,b)+\frac{b_1-a}{b-a}\Height_f(a,b,c)+\cost_f(a,b_1)\\
        &+\cost_f(b_1,b)+\cost_f(b,c) +\frac12(b-b_1)(c-b)\norm{\left(\derder[x]f\right)^-}_{L^\infty}.
    \end{align*}
    Observing that 
    \begin{itemize}
        \item $\Height_f(a,b_1,b)+\cost_f(a,b_1)+\cost_f(b_1,b)\leqslant\cost_f(a,b)$  by Prop.~\ref{prop:superadd},
        \item $\Height_f(a,b,c)\leqslant\frac12(b-a)(c-b)\norm{\left(\derder[x]f\right)^-}_{L^\infty}$ by Prop.~\ref{prop:Hbounds},
    \end{itemize}
    concludes the proof.    
    \end{proof}
 
  \subsection{Persistence property}
\label{ss:cost-persistency}

We show now that $\cost_f(a,c)$ is actually concentrated only on the countable family of maximal open intervals in $(a,c)$ where $f> \conv_{[a,c]}f$. This property will ensure that the functional $\Upsilon$   is constant along a solution to a Riemann problem and it will consequently imply also Theorem~\ref{thm:const-if-no-interaction}.

\begin{figure}[h]
			\centering
			\includegraphics{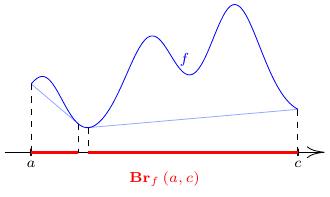}
		\caption{Definition of $\bridges_f(a,c)$.}
		\label{fig:Bridges}
	\end{figure}

\begin{definition}\label{def:bridges}
    Let $f\in C([a,c])$. We define the family of \emph{bridges of $f$ above $[a,c]$} denoted by $\bridges_f(a,c)$ as the countable family of maximal open intervals on which $f>\conv_{[a,c]}f$, see Fig.~\ref{fig:Bridges}.
\end{definition}

\begin{prop}\label{prop:persistence}
    Let $a \leqslant c$ and $f\in C^{1,1}([a,c])$. Then
    \begin{align*}
        \cost_f(a,c) = \sum_{(p,q)\in\bridges(a,c)}\cost_f(p,q).
    \end{align*}
\end{prop}
\begin{proof}
   The inequality ``$\geqslant$'' follows from superadditivity of $\cost$ (Prop.~\ref{prop:superadd}).  Let us now prove ``$\leqslant$''. 
   For every $p<q$, define the function $\pi_{p,q}(x) = p\ind{(-\infty, p]}(x)+ x\ind{(p,q)}(x)+q\ind{[q,+\infty)}(x)  $, $x\in\R$.
   A direct check shows that for any $(p,q)\in\bridges_f(a,c)$, $x,y,z\in [a,c]$, $x\leqslant y\leqslant z$, it follows
   \begin{align}
   \label{eq:bridge-singleton}
       \Height_f(\pi_{p,q}(x),\pi_{p,q}(y), \pi_{p,q}(z)) = \begin{cases}
           \Height_f(x,y,z),& y\in {(p,q)}\\
           0,& y\notin {(p,q)}.
       \end{cases}
   \end{align}
   Fix now $x,y,z\in [a,c]$, $x\leqslant y\leqslant z$. Either $y$ does not belong to any bridge, and then $\Height_f(x,y,z)=0$, or it belongs to some unique bridge. In both cases, using~\eqref{eq:bridge-singleton}, we get
   \begin{align}
   \label{prop:persistence:eq1}
       \Height_f(x,y,z) = \sum_{{(p,q)}\in\bridges_f(a,c)}\Height_f(\pi_{p,q}(x), \pi_{p,q}(y), \pi_{p,q}(z)).
   \end{align}
   
   Pick now a sequence of mergers $(b_i)_{i=1}^m  \subset [a,c]$. For each $i$, define $\bleft_i, \bright_i$ as in \eqref{def:mergecost:eq2}. For any $(p,q) \in \bridges_f(a,c)$, consider the new sequence of mergers $(\pi_{p,q}(b_i))_{i=1}^m \subset [p,q]$ and similarly define $(\pi_{p,q}(b_i))^{\rm left}, (\pi_{p,q}(b_i))^{\rm right}$. Observe that 
   \begin{align}
   \label{prop:persistence:eq2}
       \pi_{p,q}(\bleft_i) & = (\pi_{p,q}(b_i))^{\rm left},  
       & \pi_{p,q}(\bright_i) & = (\pi_{p,q}(b_i))^{\rm right}.  
   \end{align}

   Combining  the definition of merge cost in  \eqref{eq:mcostdef} with \eqref{prop:persistence:eq1}, and using \eqref{prop:persistence:eq2}, we get
   \begin{align*}
       \cost_f([a,c], (b_i)_{i=1}^m) = 
       \hspace{2mm}
       \sum_{\mathclap{(p,q)\in\bridges_f(a,c)}}
       \hspace{2mm}
       \cost_f([p,q], (\pi_{p,q}(b_i))_{i=1}^m) \leqslant 
       \hspace{2mm}
       \sum_{\mathclap{(p,q)\in\bridges_f(a,c)}}
       \hspace{2mm}
       \cost_f(p,q),
   \end{align*}
   from which the inequality ``$\leqslant$'' in the statement of the Prop. follows.   
\end{proof}
 
\subsection{Stability estimates w.r.t. second derivative of flux}
\label{ss:cost-stability}

Finally we prove some stability estimates for $\cost_f(a,c)$. In particular we show stability w.r.t. the function $f$, Prop.~\ref{prop:stability_wrt_flux}, stability w.r.t. the interval $(a,c)$, Prop.~\ref{prop:cost-lip-est},  and ``combined'' stability, both w.r.t. $f$ and w.r.t. $(a,c)$, Prop.~\ref{prop:stability}.  These properties are not used in the present paper, but they show that the definition of $\cost$ is quite robust and, as a consequence,  $\cost$ can be used in the definition of $\Upsilon$ also in the case of systems, as we  do in \cite{BaitMod2}.

First, we establish the following ``triangle inequality'' for $\cost$.

\begin{lemma}\label{lemma:costtriangle}
		If $a\leqslant c$ and $f,g\in C^{1,1}([a,c])$, then
		$
			\cost_{f+g}(a,c)\leqslant\cost_f(a,c)+\cost_g(a,c).
		$
	\end{lemma}
	\begin{proof}
		For any segment $I$, $\conv_{I}(f+g)\geqslant\conv_I f+\conv_I g$. Hence, for any $x\leqslant y\leqslant z$, $x,y,z\in[a,c]$,
		$
			\Height_{f+g}(x,y,z)\leqslant\Height_f(x,y,z)+\Height_g (x,y,z).
		$
        Applying the definition of $\cost$ concludes the proof.
	\end{proof}
    
Now we derive the following stability estimate for $\cost$ with respect to the flux function perturbation.
\begin{prop}\label{prop:stability_wrt_flux}
    If $a<c$ and $f,g\in C^{1,1}([a,c])$, then
    \begin{align}
        \Abs{\cost_f(a,c)-\cost_g(a,c)}\leqslant\frac12(c-a)\norm{\derder[x] f - \derder[x] g}_{L^1}.\label{prop:stability_wrt_flux_eq1}
    \end{align}
\end{prop}
\begin{proof}
    From Lemma~\ref{lemma:costtriangle} we get $
        - \cost_{g-f}(a,c) \leqslant \cost_f(a,c) - \cost_g(a,c) \leqslant \cost_{f-g}(a,c).$
    Then apply Prop.~\ref{prop:cost_Lone_est} to conclude to $g-f$ and $f-g$. 
\end{proof}

Now we provide the stability estimate with respect to the change of the interval.
We highlight that this estimate is motivated by the opposite-sign interactions and is fundamentally connected to Prop.~\ref{prop:subadd}.
\begin{prop}\label{prop:cost-lip-est}
    Let $a\leqslant b\leqslant c$, $f\in C^{1,1}([a,c])$. Then 
    \begin{gather*}
        \cost_f(a,b)\leqslant\cost_f(a,c)\leqslant\cost_f(a,b)+\frac12\norm{\left(\derder f\right)^-}_{L^
    \infty}\left(\frac{b+c}2-a\right)(c-b).
    \end{gather*}
\end{prop}
\begin{proof}
    The left inequality is obvious from the superadditive estimate (Prop.~\ref{prop:superadd}). For the right one, using the sub-additive estimate (Prop.~\ref{prop:subadd}) and the quadratic bound (Prop.~\ref{prop:cost_Linfty_est}), we get
    \begin{gather*}
        \cost_f(a,c)\leqslant\cost_f(a,b) + \cost_f(b,c) +
        \frac12\norm{\left(\derder f\right)^-}_{L^
    \infty} (c-b)(b-a)\leqslant\\
        \leqslant\cost_f(a,b) +\norm{\left(\derder f\right)^-}_{L^
    \infty}\left(\frac14(c-b)^2+\frac12(c-b)(b-a)\right).\qedhere
    \end{gather*}
\end{proof}
Finally, we state a stability estimate with respect to simultaneous changes of flux and the interval. This is a natural setting for the system case.
\begin{prop}\label{prop:stability}
    Let $a,b_1,b_2\in\BR$. Denote $I_i =[\min\{a,b_i\},\max\{a,b_i\}] $, $i = 1,2$. Suppose that $f_i\in C^{1,1}(I_i)$, $i = 1,2$. Then
    \begin{gather*}
        \Abs{\cost_{f_1}(a,b_1) - \cost_{f_2}(a,b_2)}\leqslant\\
        \leqslant\frac12\abs{I_1\cap I_2}\cdot\norm{\frac{\partial ^2f_1}{\partial\tau^2}-\frac{\partial ^2f_2}{\partial\tau^2}}_{L^1(I_1\cap I_2)} +\\
        +\frac12\max\left\{\norm{\frac{\partial ^2f_1}{\partial\tau^2}}_{L^\infty},\norm{\frac{\partial ^2f_2}{\partial \tau^2}}_{L^\infty}\right\}\cdot\frac{\abs{I_1}+\abs{I_2}}2\cdot \abs{b_2-b_1}.
    \end{gather*}
\end{prop}
\begin{proof}
If $a = b_1$ or $a = b_2$, the statement follows from Prop.~\ref{prop:cost_Linfty_est}. If $b_1 = b_2$, it reduces to Prop.~\ref{prop:stability_wrt_flux}.

    If $a,b_1,b_2$ are distinct, without loss of generality it is enough to consider two cases:  $a<b_1<b_2$ and $b_1<a<b_2$.

    \textit{Case $a<b_1<b_2$.}
    First, we have $I_1\cap I_2 = I_1$ and
    \begin{gather*}
        \Abs{\cost_{f_1}(a,b_1)-\cost_{f_2}(a,b_1)}\leqslant \frac12\abs{I_1\cap I_2}\cdot\norm{\frac{\partial ^2f_1}{\partial\tau^2}-\frac{\partial ^2f_2}{\partial\tau^2}}_{L^1(I_1\cap I_2)} 
    \end{gather*}
    by Prop.~\ref{prop:stability_wrt_flux}.
    Then, by Prop.~\ref{prop:cost-lip-est},
    \begin{gather*}
        \Abs{\cost_{f_2}(a,b_1)-\cost_{f_2}(a,b_2)}\leqslant \frac12\norm{\frac{\partial ^2f_2}{\partial \tau^2}}_{L^\infty}\cdot\frac{\abs{I_1}+\abs{I_2}}2\cdot(b_2-b_1).
    \end{gather*}
    
    \textit{Case $b_1<a<b_2$.}
    By Prop.~\ref{prop:cost_Linfty_est}, we have
    \begin{gather*}
        \Abs{\cost_{f_1}(a,b_1)-\cost_{f_2}(a,b_2)}\leqslant\max\{\cost_{f_1}(a,b_1),\cost_{f_2}(a,b_2)\}\leqslant\\
        \leqslant \frac14\max\left\{\norm{\frac{\partial ^2f_1}{\partial\tau^2}}_{L^\infty},\norm{\frac{\partial ^2f_2}{\partial \tau^2}}_{L^\infty}\right\} \left(\abs{I_1}+\abs{I_2}\right)^2.\qedhere
    \end{gather*}

\end{proof}

\subsection{L.s.c. w.r.t. uniform convergence of fluxes}
\label{ss:lsc-sup}

The following result, which is interesting in itself, will be used only in Section~\ref{s:semi-explicit}, namely in the proof of Theorem~\ref{th:cost-stress-new}. 

\begin{lemma}
    \label{lemma:lsc-sup}
    Let $f_n, f \in C([a,c])$ with $f_n \to f$ uniformly. Then
    \begin{equation*}
        \cost_f([a,c]) \leqslant \liminf_{n \to \infty} \cost_{f_n} ([a,c]).
    \end{equation*}
\end{lemma}
\begin{proof}
    By definition, $\Height_{f_n}(p,b,q) \to \Height_f(p,b,q)$, for all $a \leqslant p \leqslant b \leqslant q \leqslant c$. Therefore $\cost_{f_n} ([a,c], (b_j)_j) \to \cost_f([a,c], (b_j)_j)$, from which the conclusion follows. 
\end{proof}

\section{Nonnegativity, quadratic bound and interaction estimate for scalar conservation laws}
\label{s:scalar}
 
We now prove Theorem~\ref{thm:quadr-bd} in Section~\ref{ss:scalar-quadratic-bound} and Theorem~\ref{thm:interaction-estimate} in Sections~\ref{ss:scalar-same-sign} and~\ref{ss:scalar-opposite-sign}. More precisely, we prove Theorem 
\ref{thm:interaction-estimate} only in the two cases $u^L \leqslant u^M \leqslant u^R$ and $u^L \leqslant u^R \leqslant u^M$, corresponding to a same-sign interaction (Section~\ref{ss:scalar-same-sign}) and to an opposite-sign interaction (Section~\ref{ss:scalar-opposite-sign}). The other four cases (see \eqref{eq:def-delta-sigma-scalar-2}) can be treated in a completely analogous way.

\subsection{Nonnegativity and quadratic bound}
\label{ss:scalar-quadratic-bound}

\begin{proof}[Proof of Theorem~\ref{thm:quadr-bd} in the scalar case]
    By the non-negativity of $\cost$, we immediately get $\Upsilon(u)\leqslant\frac14\norm{\derder f}_{L^\infty}\TV(u)^2$. Then, applying Prop.~\ref{prop:cost_Linfty_est} to each $\cost$ term in~\eqref{UpsilonScalar}, we obtain
    \begin{align*}
        \Upsilon(u) &\geqslant \norm{\derder f}_{L^\infty} \Bigg[\frac14\TV(u)^2-\frac14\sum_{\substack{ |u_x|(\{x\})>0}}  (u_-(x)-u_{+}(x))^2\Bigg]\\
        &\geqslant\norm{\derder f}_{L^\infty} \Bigg[\frac14\TV(u)^2-\frac14\TV(u)\sum_{\substack{ |u_x|(\{x\})>0}}  \abs{u_-(x)-u_{+}(x)}\Bigg]\geqslant 0.\qedhere
    \end{align*}
\end{proof}

\subsection{Same sign interactions}
\label{ss:scalar-same-sign}
\begin{proof}[Proof of Theorem~\ref{thm:interaction-estimate} in the scalar case, $u^L\leqslant u^M\leqslant u^R$.] We adopt the notation\linebreak $a = u^L, b= u^R, c= u^M$. Observe that $\TV(u) = \TV(\hat u)$. Then
\begin{align}
\label{eq:SSI-delta-upsilon}
    \Upsilon(u)-\Upsilon(\hat u) 
    & = \cost_f(a,c) - \cost_f(a,b)-\cost_f(b,c) 
    \geqslant\Height_f(a,b,c) \\
    & =\frac12 \Delta\sigma(a,b,c).
\label{eq:SSI-delta-sigma}
\end{align}
where the inequality is due to Prop.~\ref{prop:superadd}.    
\end{proof}

\subsection{Opposite sign interactions}
\label{ss:scalar-opposite-sign}

\begin{proof}[Proof of Theorem~\ref{thm:interaction-estimate} in the scalar case, $u^L\leqslant u^R\leqslant u^M$.]
Let us first discuss the\linebreak idea of the proof. In the opposite sign interaction, the (square of the) total variation decreases, and its decrease is in fact enough to control $\Delta\sigma$ alone. However, we have the negative $\cost$ part in the functional $\Upsilon$ that accumulates, in a sense, the change of speed due to the same-sign interactions, and this part \emph{increases} across the opposite sign interaction. The non-trivial property of $\cost$ (coming from Prop.~\ref{prop:subadd}) is that this increase is also dominated by the decrease of total variation squared. Thus, we exploit twice the decrease of total variation, see~\eqref{eq:OSI-delta-sigma} and~\eqref{eq:OSI-decay-cost} below.

Differently than in the same sign interaction, we adopt here the notation $a = u^L, b= u^R, c= u^M$.
First, notice that $\TV(u)+\TV(\hat u)\geqslant 2(c-a)$, $\TV(u)-\TV(\hat u) = 2(c-b)$, so that
\begin{align}
    \TV(u)^2-\TV(\hat u)^2 \geqslant 4(c-a)(c-b)\geqslant 4(b-a)(c-b)\label{eq:OSI-decay-TV-sq}.
\end{align}
By Prop.~\ref{prop:Hbounds} and~\eqref{eq:OSI-decay-TV-sq}, we get 
\begin{align}
     \Delta\sigma(a,c,b) & = \Height_f(a,c,b) \leqslant\frac12\norm{\derder f}_{\mathclap{\hspace{4mm}L^\infty}}(c-b)(b-a) \nonumber \\ & \stackrel{\eqref{eq:OSI-decay-TV-sq}}  
      \leqslant\frac18\norm{\derder f}_{\mathclap{\hspace{4mm}L^\infty}}(\TV(u)^2-\TV(\hat u)^2 )
    .\label{eq:OSI-delta-sigma}
\end{align}
Next, applying Prop.~\ref{prop:subadd} to $\cost_f(a,c)$ and Prop.~\ref{prop:cost_Linfty_est} to $\cost_f(b,c)$ and $\cost_f(c,b)$, we get a lower bound for the change of $\cost$ part of $\Upsilon$: 
\begin{align}
    &\cost_f(a,b)-\cost_f(a,c)-\cost_f(c,b)\nonumber\\
    &\geqslant-\cost_f(b,c)-\cost_f(c,b) -\frac12\norm{\derder f}_{\mathclap{\hspace{4mm}L^\infty}}(c-b)(b-a)\nonumber \\
    &\geqslant - \frac12\norm{\derder f}_{\mathclap{\hspace{4mm}L^\infty}}(c-b)(c-a)
    \stackrel{\eqref{eq:OSI-decay-TV-sq}}\geqslant - \frac18\norm{\derder f}_{\mathclap{\hspace{4mm}L^\infty}}(\TV(u)^2-\TV(\hat u)^2 )
    .
    \label{eq:OSI-decay-cost}
\end{align}
Finally, we obtain
\begin{align*}
    \Upsilon(u)-\Upsilon(\hat u)& = \overbrace{\frac18\norm{\derder f}_{\mathclap{\hspace{4mm}L^\infty}}(\TV(u)^2-\TV(\hat u)^2 )}^{\smash{\stackrel{\eqref{eq:OSI-delta-sigma}}\geqslant\Delta\sigma(a,c,b)}}\\
    +\vphantom{\norm{\derder f}}\raisebox{-.5cm}{~}&\smash{\underbrace{\frac18\norm{\derder f}_{\mathclap{\hspace{4mm}L^\infty}}(\TV(u)^2-\TV(\hat u)^2 )+\cost_f(a,b)-\cost_f(a,c)-\cost_f(c,b)}_{\stackrel{\eqref{eq:OSI-decay-cost}}\geqslant 0}}
\end{align*}
that concludes the proof.
\end{proof}

\section{\texorpdfstring{$\Upsilon$}{Upsilon} is constant if no interaction takes place}
\label{s:const-no-int}

In this section we prove Theorem~\ref{thm:const-if-no-interaction}. The proof is actually a consequence of the following proposition.  

\begin{prop}\label{prop:scalar-RP-consistency}
     Consider the Riemann problem with initial datum
    \begin{align*}
        \bar u(x) = \begin{cases}
            a, &x<0\\
            c, &x\geqslant 0.
        \end{cases}
    \end{align*}
    Then for all $t>0$ we have $\Qcost(S_t \bar u)) = \Qcost(\bar u)$.
\end{prop}
\begin{proof}
    W.l.o.g., assume $a<c$. Let $t>0$, denote $v= S_t \bar u$.
    
    The function $\bar u$ has the only jump point $x = 0$. By Prop.~\ref{prop:persistence}, we have
    \begin{align}
        \Qcost(\bar u) = \cost_f({\bar u}_-(0), {\bar u}_+(0)) = \cost_f(a,c) = \sum_{(p,q)\in \bridges_f(a,c)}\cost_f(p,q).\label{prop:scalar-RP-consistency:eq1}
    \end{align}
    Denote now $\mathcal{S}_\lambda=\left\{
        (p,q)\in\bridges_f(a,c)\colon\frac{f(q)-f(p)}{q-p} = \lambda
        \right\}$, $\lambda\in\R$. Notice that 
        $\mathcal{S}_\lambda=\varnothing$ for all but countably many $\lambda$ and
        \begin{equation}
            \label{prop:scalar-RP-consistency:eq3}
            \bridges_f(a,c) = \bigsqcup_{\lambda} \mathcal S_\lambda.
        \end{equation}
    By definition of the solution to the Riemann problem, for each jump point $x$ of $v$, we have $\bridges_f(v_-(x), v_+(x)) = \mathcal{S}_{x/t}$ and hence again by Prop.~\ref{prop:persistence}
    \begin{align}
        \cost_f(v_-(x), v_+(x)) = \sum_{(p,q)\in \mathcal{S}_{x/t}}\cost_f(p,q)\label{prop:scalar-RP-consistency:eq2},
    \end{align}
    while, for continuity points $x$ of $v$, $\cost_f(v_-(x), v_+(x)) = 0$ and $\mathcal{S}_{x/t} = \varnothing$, so that~\eqref{prop:scalar-RP-consistency:eq2} holds as well. 
    Summing~\eqref{prop:scalar-RP-consistency:eq2} over all $x$, using~\eqref{prop:scalar-RP-consistency:eq3} and comparing with~\eqref{prop:scalar-RP-consistency:eq1}, we conclude.
\end{proof}

\begin{proof}[Proof of Theorem~\ref{thm:const-if-no-interaction}]
    Consider two adjacent jumps of $u$. Any speed originating from the first jump differs from any speed originating from the second by at most
    \begin{equation*}
            \osc\limits_{[\inf u, \sup u]} \left( \der f \right).
    \end{equation*} 
    Also, the jumps are separated spatially by at least $L_0$.
    Thus, before time $t_0$, wavefronts starting (at time zero) from different points can not interact.
    Therefore,
    for all $t<t_0$, $\TV(S_t \bar u) = \TV(\bar u)$ and, by Prop.~\ref{prop:scalar-RP-consistency}, also $\Qcost(S_t \bar u) = \Qcost(\bar u)$. 
\end{proof}

\appendix

\section{A semi-explicit integral formula for \texorpdfstring{$\cost$}{cost}}
\label{s:semi-explicit}

We defined, in Section~\ref{ss:cost-def}, $\Height_f(a,b,c)$, $\cost_f([a,c], (b_j)_j)$, $\cost_f(a,c)$ for $f \in C^{1,1}([a,c])$. Notice that the very same definitions (but not every property proven in Section~\ref{s:cost}) make sense under the weaker assumption $f \in C([a,c])$. 

We now prove a formula (Theorem~\ref{th:cost-stress-new}) that expresses the quantity $\cost_f(a,c)$, $f\in C([a,c])$ (defined in Def.~\ref{def:cost} as a supremum) in terms of a deterministic procedure, which provides a constructive way to compute $\cost$. This formula is not used anywhere in the paper,\footnote{Except in Remark~\ref{rmk:prec-consts} about optimality of constants in the bounds for $\cost$.} but we think it is interesting in itself.

\smallskip

Let thus $f \in C([a,c])$. We define the auxiliary quantity 
\begin{align}\label{eq:stress-def}
    \stress_{[a,c]}(f) := \sum_{I\in\bridges_f(a,c)}\abs{I}^2.
\end{align}
Observe that $\stress_{[a,c]}(f) \leqslant (c-a)^2$. We have the following

\begin{theorem}\label{th:cost-stress-new}
    Let $f\in C([a,c])$. Then
    \begin{align}
    \label{th:cost-stress-new-eq0}
        \cost_f(a,c) = \frac12\int_0^{+\infty}\stress_{[a,c]}(f(x)+\kappa x^2)\,d\kappa.
    \end{align}
\end{theorem}

Observe that the integrand above is non-increasing in $\kappa$. For $f \in C^{1,1}([a,c])$, if  $\kappa \geqslant \frac12 \norm{\derder f }_{L^\infty}$, the integrand is zero.\footnote{This provides an alternative proof of Prop.~\ref{prop:cost_Linfty_est}.}  For a generic continuous $f$, the integrand may be nonzero for all $\kappa \in [0,+\infty)$. 

\begin{figure}[h]
    \centering
    \includegraphics{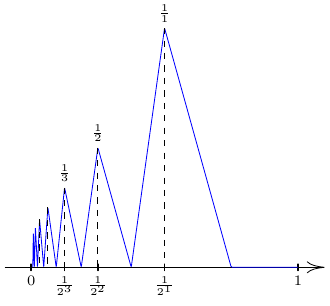}
    \caption{To Remark~\ref{rmk:ctrex}.}
    \label{fig:ctrex}
\end{figure}

\begin{remark}\label{rmk:ctrex}
    Under the current assumption $f \in C([a,c])$, 
    the $\cost$ quantity may now be infinite, in contrast to Prop.~\ref{prop:cost_Lone_est} and Prop.~\ref{prop:cost_Linfty_est}, where it was assumed $f \in C^{1,1}([a,c])$. For example, consider the function~$f$ in Fig.~\ref{fig:ctrex}. It holds $f \in C([0,1])$ and $\cost_{[0,1]} f$ is bounded from below by the sum of the length of the dashed segments $\sum_n 1/n = + \infty$. Theorem~\ref{th:cost-stress-new} takes this possibility into account as well. 
\end{remark}

The strategy of the proof is to begin with a study of the simplified discrete problem (Section~\ref{ss:PL-appr}) in which the function $f$ is piecewise linear and the merger cost w.r.t. its knots is maximized over all knots orderings. Then we proceed by passing to a limit w.r.t. the fineness of the knot set (Section~\ref{ss:cont-case}).

\subsection{The piecewise linear approximations}
\label{ss:PL-appr}
For $X = \{x_1<\ldots<x_m\}$, denote by $\PL(X)$ the set of all piecewise linear functions $\varphi\colon[x_1,x_m]\to\R$ with knots $\{x_i\}_i$, i.e. continuous functions whose graphs consist of straight line segments with ends on the lines $x = x_j$.
We define now the discrete analog of~\eqref{eq:stress-def}.
\begin{definition}\label{def:discr-stress} Let $X = \{x_1<\ldots<x_m\}$ and $\varphi\in\PL(X)$. Put
    \begin{align*}
        \stress_{X}(\varphi) = \sum_{\substack{(x_i, x_{i+s})\in\\ \bridges_\varphi(x_1,x_m)}}\Big[(x_{i+s}-x_i)^2-\sum_{j=1}^{s}(x_{i+j} - x_{i+j-1})^2\Big].
    \end{align*}
\end{definition}

\noindent Let us also introduce  the piecewise linear interpolation $\discrpar_{X}$ of the function $x \mapsto x^2$, i.e.
\begin{align*}
\discrpar_{X}\in\PL(X),\quad
\discrpar_{X}(x_j) = x_j^2,\quad \text{for } X = \{x_1 < \dots < x_m\}.
\end{align*}

The goal of this section is to prove the following result. 
\begin{prop}\label{cor:optcost}
     Let $X =\{x_1<\ldots<x_m\}$, $\varphi\in\PL(X)$. Let $(b_{j})_{j=1}^{m-2}$ be a permutation of $(x_i)_{i=2}^{m-1}$. Then
    \begin{align*}
    \cost_{\varphi}([x_1,x_m], (b_j)_{j=1}^{m-2})\leqslant \frac12 \int_0^{+\infty}\stress_{X}(\varphi+\kappa\discrpar_X)\,d\kappa.
    \end{align*}
    Moreover, there exists a permutation for which the equality holds. (The optimal permutation is given by the greedy algorithm from Lemma~\ref{lemma:greedHest} below.)
\end{prop}
\begin{remark}
    One observes that Prop.~\ref{cor:optcost} is reminiscent to Theorem~\ref{th:cost-stress-new} while def.~\ref{def:discr-stress} is slightly more technical than~\eqref{eq:stress-def}. In particular, the squares of elementary segments inside bridges are subtracted. One may informally interpret this as the exclusion of the inaccessible regions: in the discrete setting of Prop.~\ref{cor:optcost}, we are only allowed to take mergers from $X$.
\end{remark}

First, we prove a simple auxiliary lemma.
\begin{lemma}\label{lemma:bridges_splitting}
    Let $X\subset\R$ be finite and $\varphi\in\PL(X)$. Let $a<b<c$ be any three elements of~$X$ and $b\in (p,q)\in\bridges_\varphi(a,c)$. Then
    \begin{align*}
        (q-b)(b-p)\leqslant 
        \frac{1}{2}\Big(\stress_{X \cap [a,c]}(\varphi)- \stress_{X\cap[a,b]}(\varphi)
        -\stress_{X\cap[b,c]}(\varphi) \Big).
    \end{align*}
\end{lemma}
\begin{proof}
    Note that the expression inside the brackets in the r.h.s. is in fact an algebraic sum of squares of several interval lengths. After all cancellations, it reduces to the form
    \begin{align*}
        (q-p)^2 - \boxed{
        \begin{minipage}{5cm}
            sum of squared lengths of some disjoint subintervals inside $[p,b]$
        \end{minipage}
        }
        -
        \boxed{
        \begin{minipage}{5cm}
            sum of squared lengths of some disjoint subintervals inside $[b,q]$
        \end{minipage}
        }.
    \end{align*}
    The boxes in the r.h.s. are bounded from above by $(b-p)^2$ and $(q-b)^2$, respectively, that concludes the proof.
\end{proof}

\begin{lemma}\label{lemma:Hkappader}
    Let $X\subset\R$ be finite and $\varphi\in\PL(X)$. Let $a<b<c$ be any three elements of~$X$. Then for every $\kappa$, there exists the right derivative $\partial_\kappa^+ \Height_{\varphi+\kappa\discrpar_X}(a,b,c)$ and
    \begin{enumerate}[(i)]
        \item if $\Height_{\varphi+\kappa\discrpar_X}(a,b,c)=0$, then $ \partial_\kappa^+ \Height_{\varphi+\kappa\discrpar_X}(a,b,c) = 0$;\label{lemma:Hkappader:st1}
        \item if $\Height_{\varphi+\kappa\discrpar_X}(a,b,c)>0$, then $
        \partial_\kappa^+ \Height_{\varphi+\kappa\discrpar_X}(a,b,c) = -(q-b)(b-p)$,
        where the interval $(p,q)\in\bridges_{\varphi+\kappa\discrpar_X}(a,c)$ is the unique bridge such that $b\in(p,q)$.\label{lemma:Hkappader:st2}
    \end{enumerate} 
\end{lemma}
\begin{proof} 
    Notice that if $\Height_{\varphi+\kappa\discrpar_X}(a,b,c) $ vanishes for some $\kappa = \kappa_0$, it also vanishes for every $\kappa>\kappa_0$. This proves~(\ref{lemma:Hkappader:st1}). Let us prove~(\ref{lemma:Hkappader:st2}).
    Since $\varphi,\discrpar_X\in\PL(X)$, the interval $(p,q)$ remains a bridge for $\varphi + (\kappa + \delta) \discrpar_X$,  i.e. $b \in (p,q)\in\bridges_{\varphi+(\kappa+\delta)\discrpar_X}(a,c)$, for $\delta>0$ small enough. Hence, for $f \in\{\varphi + \kappa \discrpar_X, \varphi + (\kappa + \delta) \discrpar_X \} $ it holds
    \begin{align*}
        \Height_f(a,b,c) = \Height_f(p,b,q) = f(b) - \tfrac{q-b}{q-p}f(p)-\tfrac{b-p}{q-p}f(q),\quad 
    \end{align*}
    and so
    \begin{align*}
        &\Height_{\varphi+(\kappa+\delta)\discrpar_X}(a,b,c) - \Height_{\varphi+\kappa\discrpar_X}(a,b,c) 
        =\delta{\smash{\Big[\discrpar_X(b) - \tfrac{q-b}{q-p}\, \discrpar_X(p)-\tfrac{b-p}{q-p}\, \discrpar_X(q)\Big]}}\\
        &=-\delta(q-b)(b-p),
    \end{align*}
    where the last equality holds because $p,b,q \in X$. 
\end{proof}
\begin{lemma}\label{cor:Hkappader}
    In the setting of Lemma~\ref{lemma:Hkappader},
    \begin{align*}
        0&\geqslant \partial_\kappa^+ \Height_{\varphi+\kappa\discrpar_X}(a,b,c)\\
        &\geqslant \frac{1}{2}\Big(\stress_{X\cap[a,b]}(\varphi+\kappa\discrpar_X)+\stress_{X\cap[b,c]}(\varphi+\kappa\discrpar_X)-\stress_{X \cap [a,c]}(\varphi+\kappa\discrpar_X) \Big).
    \end{align*}
\end{lemma}
\begin{proof}
    Combine Lemma~\ref{lemma:Hkappader} with Lemma~\ref{lemma:bridges_splitting}.
\end{proof}
\begin{lemma}\label{lemma:coststressest}
    Let $X =\{x_1<\ldots<x_m\}$, $\varphi\in\PL(X)$. Let $(b_{j})_{j=1}^{m-2}$ be an arbitrary permutation of $(x_i)_{i=2}^{m-1}$.  Then there exists the right derivative
    \begin{align}         \partial_{\kappa}^+\cost_{\varphi+\kappa\discrpar_X}([x_1,x_m], (b_j)_{j=1}^{m-2})\geqslant -\frac12 \stress_{X}(\varphi+\kappa\discrpar_X).
         \label{lemma:coststressest:eq1}
    \end{align}
\end{lemma}
\begin{proof}
        Using the formula~\eqref{eq:mcostdef} and Lemma~\ref{cor:Hkappader}, we get
        \begin{align*}
            &\partial_\kappa^+\cost_{\varphi+\kappa\discrpar_X}([x_1,x_m], (b_j)_{j=1}^{m-2}) =  \sum_{j=1}^{m-2}\partial_\kappa^+\Height_{\varphi+\kappa\discrpar_X}(\bleft_j, b_j, \bright_j)\\
            & \geqslant\frac12 \sum_{j=1}^{m-2}\Big(\stress_{{X\cap[\bleft_j,b_j]}}+\stress_{{X\cap[b_j,\bright_j]}}-\stress_{{X \cap [\bleft_j,\bright_j]}}\Big)(\varphi+\kappa\discrpar_X),
        \end{align*}
        where $\bleft_j$, $\bright_j$ are as in Def.~\ref{def:mergecost}.
        By the telescopic summation trick, analogously to the computation in the proof of Lemma~\ref{lemma:psi_triangle}, one reduces the latter sum to
        \begin{align*}
        \Big[\stress_{\{x_1,x_2\}}+\ldots+\stress_{\{x_{m-1},x_m\}}-\stress_{X}\Big](\varphi+\kappa\discrpar_X), 
        \end{align*}
        where all the terms but the last one vanish. This concludes the proof.
\end{proof}
In Lemmas~\ref{lemma:greedHest} and~\ref{lemma:optdcost} below, we show that some ordering fulfills the inequalities in Lemmas~\ref{cor:Hkappader} and \ref{lemma:coststressest} respectively.
\begin{lemma}\label{lemma:greedHest}
    Let $X = \{x_1<\ldots<x_m\}$, $\varphi\in\PL(X)$. Let the permutation $(b_j)_{j=1}^{m-2}$ of $(x_{i})_{i=2}^{m-1}$ be obtained by the following \emph{greedy algorithm}:
    \begin{itemize}
        \item we choose inductively $b_1,\ldots, b_{m-2}$;
        \item on $j^{\text{th}}$ step, $b_j$ is chosen such that the quantity 
        \begin{align}
        \frac{\Height_{\varphi}\left(\bleft_j, b_j, \bright_j\right)}{(\bright_j-b_j)(b_j-\bleft_j)}\label{eq:greedy}
        \end{align}
        is minimal over all possible remaining choices, where $\bleft_j$, $\bright_j$ are as in Def.~\ref{def:mergecost}.
    \end{itemize}
    Then
    \begin{align}
&\partial_{\kappa}^+\Height_{\varphi+\kappa\discrpar_X}\left(\bleft_j, b_j, \bright_j\right)\nonumber\\
&=
\frac12\Big(\stress_{X\cap[\bleft_j,b_j]}+\stress_{X\cap[b_j,\bright_j]}-\stress_{X \cap [\bleft_j, \bright_j]}\Big)(\varphi+\kappa\discrpar_X) .\label{lemma:greedHest:eq1}
    \end{align}
\end{lemma}
\begin{proof}
     Fix $\kappa\geqslant 0$. Assume first that the r.h.s. of~\eqref{lemma:greedHest:eq1} equals $0$. Then Lemma~\ref{cor:Hkappader} implies the l.h.s. to equal $0$ as well.

    Assume now that the r.h.s. of~\eqref{lemma:greedHest:eq1} is negative.
    Then $\Height_{\varphi+\kappa\discrpar_X}\left(\bleft_j, b_j, \bright_j\right)\allowbreak{>0}$ and a fortiori $\Height_\varphi(\bleft_j, b_j, \bright_j) >0$. 
    By the ``greedy'' choice of $b_j$, for all $x\in X\cap (\bleft_j, \bright_j)$ we have   
    \begin{align*}
        \frac{\varphi(x) - \big(\conv_{[\bleft_j, \bright_j]}\varphi\big)(x) }{(\bright_j -x )(x-\bleft_j)}\geqslant\frac{\Height_{\varphi}\left(\bleft_j, b_j, \bright_j\right)}{(\bright_j-b_j)(b_j-\bleft_j)}
    \end{align*}
    i.e.
    \begin{align}
        \varphi(x) \geqslant \big(\conv_{[\bleft_j, \bright_j]}\varphi\big)(x) + \frac{\Height_{\varphi}\left(\bleft_j, b_j, \bright_j\right)}{(\bright_j-b_j)(b_j-\bleft_j)}(\bright_j -x )(x-\bleft_j),
        \label{lemma:greedHest:eq2}
    \end{align}
      with equality holding for $x = b_j$. Notice that equality holds in~\eqref{lemma:greedHest:eq2} also for $x = \bleft_j, \bright_j$. Since the second term in r.h.s. of~\eqref{lemma:greedHest:eq2} is strictly positive on the internal knots $X\cap(\bleft_j, \bright_j)$, the function $\varphi\in\PL(X)$ is strictly bigger than its convex envelope on  $X \cap (\bleft_j,\bright_j)$  and therefore, $\conv_{[\bleft_j, \bright_j]}\varphi$ is affine. The r.h.s. in~\eqref{lemma:greedHest:eq2} is thus a quadratic function, that we denote by $P(x)$,  see Fig.~\ref{fig:greedHest}(a).
     Then
    \begin{align}
        \left(\varphi+\kappa\discrpar_X\right)(x)\geqslant P(x) +\kappa x^2,  \quad x\in X\cap [\bleft_j, \bright_j].\label{lemma:greedHest:eq3}
    \end{align}
    \begin{figure}[h]
        \centering
        \begin{minipage}{0.45\textwidth}
            \centering
            \vspace{6.3mm}\includegraphics{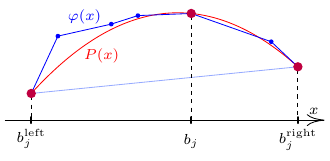}\\
            (a)
        \end{minipage} \hfill
        \begin{minipage}{0.45\textwidth}
            \centering
            \includegraphics{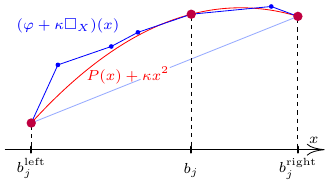}\\
            (b)
        \end{minipage}
        \caption{To the proof of Lemma~\ref{lemma:greedHest}.}
        \label{fig:greedHest}
    \end{figure}
    with equality for $x = \bleft_j, b_j, \bright_j$. Since we are in the case  $\Height_{\varphi + \kappa \discrpar_X}(\bleft_j, b_j, \bright_j)\allowbreak>0$, the parabola $P(x) + \kappa x^2$ is strictly concave, see Fig.~\ref{fig:greedHest}(b). Taking this into account, we deduce, using~\eqref{lemma:greedHest:eq3},
    \begin{align*}
    & \bridges_{\varphi+\kappa\discrpar_X}{(\bleft_j,b_j)} = \{(\bleft_j,b_j)\},\quad \bridges_{\varphi+\kappa\discrpar_X}{(b_j,\bright_j)} = \{(b_j,\bright_j)\}\\
    &\bridges_{\varphi+\kappa\discrpar_X}{(\bleft_j,\bright_j)} = \{(\bleft_j,\bright_j)\}.
    \end{align*}
    This, together with Lemma~\ref{lemma:Hkappader}, item~\eqref{lemma:Hkappader:st2}, implies~\eqref{lemma:greedHest:eq1}.    
\end{proof}
\begin{lemma}\label{lemma:optdcost}
    Let the assumptions of Lemma~\ref{lemma:greedHest} hold. Then
     \begin{align}
     \partial_{\kappa}^+\cost_{\varphi+\kappa\discrpar_X}([x_1,x_m], (b_j)_{j=1}^{m-2})= -\frac12 \stress_{X}(\varphi+\kappa\discrpar_X).
         \label{lemma:optdcost:eq1}
    \end{align}
\end{lemma}
\begin{proof}
Mimic the proof of Lemma~\ref{lemma:coststressest} with equality instead of inequality (due to Lemma~\ref{lemma:greedHest}).
\end{proof}
From the estimates we have just developed, we are now able to derive Prop.~\ref{cor:optcost}.
\begin{proof}[Proof of Prop.~\ref{cor:optcost}]
    Fix an arbitrary permutation $(b_j)_{j=1}^{m-2}$ of $(x_i)_{i=2}^{m-1}$. Notice that $\varphi+\kappa\discrpar_X$ is convex for $\kappa$ large enough, so that both functions
    \begin{align*}
        F(\kappa):= \cost_{\varphi+\kappa\discrpar_X}([x_1,x_{m}], (b_j)_{j=1}^{m-2}),\quad G(\kappa):=\stress_{X}(\varphi+\kappa\discrpar_X)
    \end{align*} are zero for $\kappa$ large enough. By construction, $F$ is continuous, and moreover, since $f,\discrpar_X\in\PL(X)$, $F$ is piecewise linear and $G$ is piecewise constant. Thus
    \begin{align*}
         \cost_{\varphi}([x_1,x_{m}], (b_j)_{j=1}^{m-2})=F(0)=-\int\limits_{0}^{+\infty}\partial^+_\kappa F(\kappa)\,d\kappa\stackrel{\text{Lemma~\ref{lemma:coststressest}}}\leqslant\frac12\int\limits_{0}^{+\infty}G(\kappa)\, d\kappa.
    \end{align*}
    In the special case when $(b_j)_j$ is given by the greedy algorithm from Lemma~\ref{lemma:greedHest}, we have ``$=$'' instead of ``$\leqslant$'' by Lemma~\ref{lemma:optdcost}.
\end{proof}

\subsection{The continuous case}
\label{ss:cont-case}
We begin with the  approximation results for the ``discrete'' and ``continuous'' $\stress$'es (recall~\eqref{eq:stress-def} and Def.~\ref{def:discr-stress}) in Lemmas~\ref{lemma:stressapprox} and~\ref{cor:stressconv} below.
\begin{lemma}\label{lemma:stressapprox}
    Let $a<c$ and $X_N = \{a=x^N_1<\ldots<x^N_{m_N} = c\}$ be a finite subset of $\R$ for every $N = 1,2\ldots$. Suppose that $\max_i(x^{N}_{i+1}-x^N_i)\to 0$ as $N\to\infty$. Let $f\in C([a,c])$ and $f_N$ be the piecewise linear interpolation of $f$ with knots $X_N$. Then for every $\delta>0$
    \begin{align}
        \liminf_{N\to\infty} \stress_{X_N} (f_N)\geqslant\stress_{[a,c]}(f)\geqslant\limsup_{N\to\infty} \stress_{X_N} (f_N+\delta\discrpar_{X_N}).\label{lemma:stressapprox:eq1}
    \end{align}
\end{lemma}
\begin{proof}
    Notice that,
    by Def.~\ref{def:discr-stress} and~\eqref{eq:stress-def}, $\max_i(x^{N}_{i+1}-x^N_i)\to 0$ implies
    \begin{align*}
        \stress_{X_N} (f_N) - \stress_{[a,c]}(f_N) \to 0,\quad \stress_{X_N} (f_N+\delta\discrpar_{X_N}) - \stress_{[a,c]} (f_N+\delta\discrpar_{X_N})\to 0
    \end{align*}
    as $N \to \infty$,
    so we will show~\eqref{lemma:stressapprox:eq1} with discrete stresses $\stress_{X_N}$ substituted by  $\stress_{[a,c]}$.
    
    For every $h\in C([a,c])$, denote
    \begin{align*}
        d_h(x) = \dist\Big(x;[a,c]\setminus \bigcup_{I\in\bridges_h(a,c)} I\Big).
    \end{align*}
    Notice that
    \begin{subequations}
        \begin{align}
        d_h(x)&= \dist\left(x, \left\{y\colon h(y) = \big(\conv_{[a,c]}h\big)(y)\right\}\right)\label{lemma:stressapprox:eq1.51} \\
        &= \dist(x, \{y\colon d_h(y) = 0\})
        \label{lemma:stressapprox:eq1.52}
    \end{align}
    \end{subequations}
    and
    \begin{align}
        \stress_{[a,c]} (h) = 4 \int_a^c d_h(x)\, dx.\label{lemma:stressapprox:eq1.6}
    \end{align}
    We will show 
    \begin{align}
            \forall x\in [a,c]\colon\quad\liminf_{N\to\infty}d_{f_N}(x) \geqslant d_f(x)\geqslant \limsup_{N\to\infty} d_{f_N+\delta\discrpar_{X_N}}(x).\label{lemma:stressapprox:eq4}
        \end{align}
        Then~\eqref{lemma:stressapprox:eq1.6} and Fatou's Lemma implies the first inequality in~\eqref{lemma:stressapprox:eq1}, while~\eqref{lemma:stressapprox:eq1.6} and reverse Fatou's Lemma together with the bound  $0 \leqslant d_{(\cdot)} \leqslant c-a$ implies the second inequality in~\eqref{lemma:stressapprox:eq1}.
    
    Since $f_N\to f$ uniformly, we also have $f_N-\conv_{[a,c]} f_N \to f -\conv_{[a,c]} f$ uniformly and hence, by~\eqref{lemma:stressapprox:eq1.51}, the first inequality in~\eqref{lemma:stressapprox:eq4} follows.

     We will now show that, for any $z\in[a,c] $ with $d_f(z) = 0$, there exists a sequence $z_N\to z$ such that $z_N\in X_N$ and $d_{f_N+\delta\discrpar_{X_N}}(z_N) = 0$. This  implies the second inequality in~\eqref{lemma:stressapprox:eq4}  due to~\eqref{lemma:stressapprox:eq1.52}. 

    Fix $z\in [a,c]$, $d_f(z) =0$. W.l.o.g. assume $z = 0$ and $f(z) = 0$.
    There exists a number $L\in\R$ such that $
        f(x)\geqslant Lx$ for all $x$. Put $g(x) = f(x)-Lx+\delta x^2$. 
        For every $N$, chose $z_N\in X_N $ to be any minimizer of $g$ over $ X_N$. 
        
        Since $g$ is continuous and the $X_N$ is approaching every point, $g(z_N)\to 0$. By definition $g(x) \geqslant \delta x^2$ for all $x$. Then also $\abs{z_N - z} = \abs{z_N}\leqslant\delta^{-1/2}\sqrt{ g(z_N)}\to 0$. 
        
        We prove now $d_{f_N+\delta\discrpar_{X_N}}(z_N) = 0$. Being $z_N$ a minimizer, we have for every $x\in X_N$ the inequalities
        \begin{align}
            f(x)&+\delta x^2& &\geqslant & f(z_N) &+\delta z_N^2 + L(x-z_N) \nonumber\\
            &\hspace{1.7mm}\rotatebox{90}{=}& & & &\hspace{1cm}\rotatebox{90}{=} \nonumber\\
            f_N(x)&+\delta\discrpar_{X_N}(x) & &\geqslant & f_N(z_N) &+\delta z_N^2 + L(x-z_N),\label{lemma:stressapprox:eq6}
        \end{align}
        and thus, using that $f_N, \discrpar_{X_N} \in \PL(X_N)$,  the inequality in~\eqref{lemma:stressapprox:eq6} holds for every $x \in [a,c]$, with equality for $x = z_N$. Therefore the r.h.s. in~\eqref{lemma:stressapprox:eq6} is a supporting line for the graph of $f_N + \delta \discrpar_{X_N}$, and    
        this proves $d_{f_N+\delta\discrpar_{X_N}}(z_N) = 0$.
\end{proof}
\begin{lemma}\label{cor:stressconv}
    In the same setting as Lemma~\ref{lemma:stressapprox}, for almost every $\kappa\in\R$ we have
    \begin{align*}
       \stress_{X_N} (f_N+\kappa\discrpar_{X_N})\to \stress_{[a,c]}(f(x)+\kappa x^2).
    \end{align*}
\end{lemma}
\begin{proof}
Put $f_\kappa(x) = f(x) + \kappa x^2$.
Notice that
$\kappa\mapsto \stress_{[a,c]}(f_\kappa)$ is non-increasing. Hence, for a.e. $\kappa\in \R$,
$
    \stress_{[a,c]}(f_{\kappa-\delta})\to \stress_{[a,c]}(f_\kappa)
$ follows as $\delta\searrow 0$.    
    But for every $\kappa$ and $\delta>0$, by Lemma~\ref{lemma:stressapprox}, we have
    \begin{align*}
        \stress_{[a,c]}(f_{\kappa-\delta})&\geqslant\limsup_{N\to\infty} \stress_{X_N} (f_N+\kappa\discrpar_{X_N})\\
        &\geqslant
        \liminf_{N\to\infty} \stress_{X_N} (f_N+\kappa\discrpar_{X_N})\geqslant \stress_{[a,c]}(f_\kappa).\qedhere
    \end{align*}
\end{proof}
Now we are able to derive the desired formula in Theorem~\ref{th:cost-stress-new} with auxiliary Lemma~\ref{lemma:intstressXest}.
\begin{lemma}\label{lemma:intstressXest}
    Let $X = \{a = x_1<\ldots< x_m = c\}$ and $f\in C([a,c])$. Let $f_X\in \PL(X)$ be the piecewise linear interpolation of $f$ with knots from $X$. Then
    \begin{align*}
        \frac12\int_0^{+\infty} \stress_X(f_X+\kappa\discrpar_X)\, d\kappa \leqslant\cost_f(a,c).
    \end{align*}
\end{lemma}
\begin{proof}
    Let $(b_j)_{j=1}^{m-2}$ be the permutation of $(x_i)_{i=2}^{m-1}$ given by the greedy algorithm from Lemma~\ref{lemma:greedHest}, applied to $\varphi = f_X$. Then
    \begin{align*}
        \frac12\int_0^{+\infty} \stress_X(f_X+\kappa\discrpar_X)\, d\kappa  =  \cost_{f_X}([a,c];(b_j)_{j=1}^{m-2}) \leqslant\cost_{f}([a,c];(b_j)_{j=1}^{m-2}),
    \end{align*}
    where the inequality above holds since for all $x,y,z \in X$, $\Height_{f_X}(x,y,z) \allowbreak\leqslant \Height_f(x,y,z)$. The conclusion now follows from the definition of $\cost_f(a,c)$. 
\end{proof}
\begin{proof}[Proof of Theorem~\ref{th:cost-stress-new}]
    We have the inequality
    \begin{align}
    \frac12\int_0^{+\infty} \stress_{[a,c]}(f(x)+\kappa x^2)\, d\kappa \leqslant\cost_f(a,c).\label{th:cost-stress:eq1} 
    \end{align}
    as a consequence of Lemmas~\ref{cor:stressconv},~\ref{lemma:intstressXest}, and Fatou's Lemma.

    Let us  now prove the opposite inequality. We first assume that for some $\kappa_0>0$, ${f(x)+\kappa_0x^2}$ is convex.\footnote{This is not always the case, e.g. consider $f(x) = -|x|$ on $[-1,1]$.}

    Fix $N\in\BN$. Choose finitely many points $b_1, b_2,\ldots\in [a,c]$ such that
    \begin{align}
    \label{th:cost-stress:eq5}
        \cost_f(a,c)<\cost_f([a,c], (b_j)_{j})+1/N.
    \end{align}
    We can assume w.l.o.g. that:
    \begin{enumerate}[(i)]
        \item For every $x \in [a,c]$ there is $b_i$ such that $|x - b_i| < 1/N$.
        \item For every $i$, $\conv_{[\bleft_i, \bright_i]} f$ is affine, if $\bleft_i<\bright_i$.  
    \end{enumerate}
    Indeed, if (i) is not true, one can simply enlarge the sequence $(b_j)_j$ by appending more points. This just makes the r.h.s. in~\eqref{th:cost-stress:eq5} larger. If (ii) is not true, let $p_i, q_i$ be such that $b_i\in (p_i, q_i)\in \bridges_{f}(\bleft_i,\bright_i) $. One then insert $p_i, q_i$ into the sequence $(b_j)_j$ just before $b_i$. This does not change the merger cost.

    Put $X_N = \{b_j\}_j$ and let $f_N\in\PL(X_N)$ be the piecewise linear interpolation of $f$ with knots~$X_N$. Thanks to Property (ii) above, 
    \begin{equation}
    \label{th:cost-stress:eq6}
        \Height_f(\bleft_j, b_j, \bright_j) = \Height_{f_N}(\bleft_j, b_j, \bright_j)
    \end{equation}
    and therefore also
    \begin{align*}
        \cost_f (a,c) - \frac{1}{N}\ &\stackrel{\mathclap{\eqref{th:cost-stress:eq5}}}{\leqslant}\ \cost_f([a,c], (b_j)_{j}) \stackrel{\eqref{th:cost-stress:eq6}}= \cost_{f_N}([a,c], (b_j)_{j})\\
        &\stackrel{\mathclap{\text{Prop.~\ref{cor:optcost}}}}{\leqslant }\hspace{5mm}\frac12\int_{0}^{+\infty}\stress_{X_N}(f_N+\kappa\discrpar_{X_N})\, d\kappa.
    \end{align*}
    Notice that the integrand in the last term is the piecewise linear interpolation of $f(x) + \kappa x^2$, which is convex for $\kappa \geqslant \kappa_0$. Therefore the integrand equals $0$ whenever $\kappa\geqslant\kappa_0$. We can thus use Lemma~\ref{cor:stressconv} (which is applicable because of (i) above) together with Dominated Convergence Theorem to get
    \begin{align*}
        \cost_f(a,c) \leqslant \frac12\int_0^{+\infty} \stress_{[a,c]}(f(x)+\kappa x^2)\, d\kappa 
    \end{align*}
    This concludes the proof for the case when ${f(x)+\kappa_0x^2}$ is convex.

    Consider now the general case $f\in C([a,c])$. Put
    \begin{align*}
        f_\kappa(x) = f(x) +\kappa x^2,\quad g_\kappa (x) = \Big(\conv_{[a,c]}f_{\kappa}\Big)(x)- \kappa x^2\leqslant f(x),
    \end{align*}
    see Fig.~\ref{fig:g_kappa}.
    We claim that
     \begin{align}
        \stress_{[a,c]}(g_\kappa(x)+\kappa' x^2)=\begin{cases}
             \stress_{[a,c]}(f_{\kappa'}),&\kappa'<\kappa,\\
             0,& \kappa'\geqslant\kappa.
        \end{cases}\label{th:cost-stress:eq2}
    \end{align}
    \begin{figure}[h]
        \centering
        \includegraphics{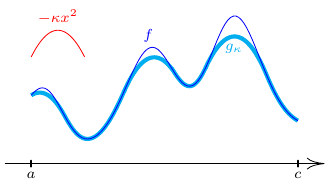}
        \caption{To the definition of $g_\kappa$.}
        \label{fig:g_kappa}
    \end{figure}
     Assume that~\eqref{th:cost-stress:eq2} is proved. Since $g_\kappa(x)+\kappa x^2$ is convex, by the first part of the proof, we can apply Theorem~\ref{th:cost-stress-new} to $g_\kappa$ and we get
     \begin{equation}
     \label{th:cost-stress:eq2.5}
         \cost_{g_\kappa}(a,c) = \frac12\int_0^{+\infty} \stress_{[a,c]}(g_\kappa(x)+\kappa' x^2)\, d\kappa'.
     \end{equation}
     We have $\norm{f - g_\kappa}_{L^\infty} =  \norm{f_\kappa-\conv_{[a,c]} f_\kappa}_{L^\infty}\allowbreak\xrightarrow{\kappa\to+\infty} 0$. Indeed, let $\rho_\e$ be a mollifier such that $p * \rho_\e =p$ for every polynomial $p$ with $\deg(p) \leqslant 2$, so that $f_\kappa * \rho_\e = (f * \rho_\e)_\kappa = f* \rho_\e + \kappa x^2  =: f_\kappa^\e$. Therefore
     \begin{align*}
         \hspace{5mm}\rule{0pt}{7mm}\smash{f_\kappa - \conv_{[a,c]} f_\kappa \leqslant 2 \| \overbrace{f_\kappa - f_\kappa^\e}^{= f - f^\e}\|_{L^\infty} + \overbrace{f_\kappa^\e - \conv_{[a,c]} f_\kappa^\e.}^{=0 \text{ if } \kappa \geqslant \frac12 \norm{\derder{f_\kappa^\e}}_{L^\infty} } }
     \end{align*}
     We use now l.s.c. of $\cost$ w.r.t. uniform convergence (Lemma~\ref{lemma:lsc-sup}) to get
    \begin{align*}        \cost_f(a,c) \leqslant\liminf_{\kappa\to+\infty}\cost_{g_\kappa}(a,c) \stackrel{\eqref{th:cost-stress:eq2.5}}= \liminf_{\kappa\to+\infty}\frac12\int_0^{+\infty} \stress_{[a,c]}(g_\kappa(x)+\kappa' x^2)\, d\kappa'.
    \end{align*}
    Using~\eqref{th:cost-stress:eq2}, we get the opposite inequality to~\eqref{th:cost-stress:eq1} and therefore conclude the proof of the theorem. 

    It remains to prove~\eqref{th:cost-stress:eq2}.
    Fix $\kappa'\geqslant 0$. If $\kappa'\geqslant \kappa$,  $g_\kappa(x) + \kappa'x^2$ is convex and thus~\eqref{th:cost-stress:eq2} follows. Assume now $\kappa'<\kappa$.
     We are going to derive the equality
     \begin{align}\label{th:cost-stress:eq3}
     \bridges_{g_{\kappa}(x)+\kappa'x^2}(a,c) = \bridges_{f_{\kappa'}}(a,c)    
     \end{align}
     from which~\eqref{th:cost-stress:eq2}  follows.

     First, we have 
     \begin{equation*}
       (\conv_{[a,c]} f_{\kappa'})(x)+(\kappa-\kappa')x^2 \leqslant f_{\kappa'}(x)+(\kappa-\kappa')x^2= f_{\kappa}(x).
     \end{equation*}
     The l.h.s. above is convex, and thus 
    \begin{align*}
        \big(\conv_{[a,c]}f_{\kappa'}\big)(x)\leqslant\big(\conv_{[a,c]}f_{\kappa}\big)(x) - (\kappa-\kappa')x^2 = \smash{\overbrace{g_\kappa(x)}^{\leqslant f(x)}}+\kappa'x^2\leqslant f_{\kappa'}(x),
    \end{align*}
    so that, taking $\conv_{[a,c]}$ on all terms, we get
    \begin{align*}
        \conv_{[a,c]}(g_\kappa(x)+\kappa' x^2) = \conv_{[a,c]}f_{\kappa'}. 
    \end{align*}
    Since $g_\kappa(x)+\kappa' x^2\leqslant f_{\kappa'}(x)$, for the justification of~\eqref{th:cost-stress:eq3}  it is sufficient now to prove the implication
    \begin{align*}
       \forall x_0\in[a,c]\colon\quad g_\kappa(x_0)+\kappa' x_0^2=\conv_{[a,c]}(f_{\kappa'})(x_0) \quad\Rightarrow\quad f_{\kappa'}(x_0) = \conv_{[a,c]}(f_{\kappa'})(x_0).
    \end{align*}
    
        Let thus $x_0$ fulfills the premise of the implication, then  there exists $L$  such that 
        \begin{align*}
            \forall  x\in[a,c]\colon\quad g_\kappa(x)+\kappa' x^2\geqslant g_\kappa(x_0)+\kappa' x_0^2+ L(x-x_0).
        \end{align*}
        Substituting the definition of $g_\kappa $ to the latter inequality, we get
        \begin{align}
        \label{th:cost-stress:eq4}
            (\conv_{[a,c]}f_{\kappa})(x)\geqslant (\conv_{[a,c]}f_{\kappa})(x_0)+(\kappa-\kappa')(x^2-x_0^2)+L(x-x_0).
        \end{align} 
        The r.h.s. in~\eqref{th:cost-stress:eq4} is an upward-opening parabola, which is less or equal than $\conv_{[a,c]}f_{\kappa}$, with equality holding  for $x=x_0$. Hence $\conv_{[a,c]} f_\kappa$ can not be affine on a neighborhood of $x_0$ and therefore
        $f_\kappa(x_0) = (\conv_{[a,c]}f_\kappa)(x_0)$.
        So 
        \begin{align*}
            f_{\kappa}(x)\geqslant(\conv_{[a,c]}f_{\kappa})(x)\smash{\stackrel{\eqref{th:cost-stress:eq4}}\geqslant} f_{\kappa}(x_0)+(\kappa-\kappa')(x^2-x_0^2)+L(x-x_0)
        \end{align*}
         and hence
    $f_{\kappa'}(x)\geqslant f_{\kappa'}(x_0)+L(x-x_0)$. We deduce $ f_{\kappa'}(x_0) = \conv_{[a,c]}(f_{\kappa'})(x_0)$.
\end{proof}
\begin{remark}
    We mentioned in Remark~\ref{rmk:quad-cubic} that there is some analogy between $\cost$ and the quantity $\mathcal{A}$ in \cite{Caravenna2008}. 
    One could also prove that if we substitute $\abs{I}^2$ by $\abs{I}^3$, then the formula in Theorem~\ref{th:cost-stress-new} gives precisely $\mathcal{A}(a,c) = \norm{f - \conv_{[a,c]} f}_{L^1([a,c])}$ instead of $\cost_{[a,c]} f$. 
\end{remark}

\bibliographystyle{plain}
\bibliography{references}

\end{document}